\documentclass[twoside,10pt,english]{smfart}

\usepackage[matrix,arrow]{xy}
\usepackage{amssymb,amsmath}
\usepackage[mathscr]{euscript}
\usepackage{mathrsfs}

\title[Local and stable homological algebra]{Local and stable homological algebra\\
in Grothendieck abelian categories}

\date{December 2006, last modified in December 2007}

\author{Denis-Charles Cisinski}
\address{LAGA\\
CNRS~(UMR 7539)\\
Institut~Galil\'ee\\
Universit\'e~Paris~13\\
\hbox{Avenue~Jean-Baptiste~Cl\'ement}\\
93430~Villetaneuse\\France}
\email{cisinski@math.univ-paris13.fr}
\urladdr{http://www.math.univ-paris13.fr/~cisinski/}

\author{Fr\'ed\'eric D\'eglise}
\address{LAGA\\
CNRS (UMR 7539)\\
Institut Galil\'ee\\
Universit\'e Paris~13\\
\hbox{Avenue~Jean-Baptiste~Cl\'ement}\\
93430~Villetaneuse\\France}
\email{deglise@math.univ-paris13.fr}
\urladdr{http://www.math.univ-paris13.fr/~deglise/}

\newtheorem{thm}{Theorem}[section]
\newtheorem{prop}[thm]{Proposition}
\newtheorem{lm}[thm]{Lemma}
\newtheorem{cor}[thm]{Corollary}

\theoremstyle{remark} 
\newtheorem{rem}[thm]{Remark}

\newtheorem{ex}[thm]{Example}

\theoremstyle{definition} 
\newtheorem{df}[thm]{Definition}
\newtheorem{num}[thm]{}
\newtheorem{paragr}[thm]{}
\numberwithin{equation}{thm}

\DeclareFontFamily{OT1}{pzc}{}
\DeclareFontShape{OT1}{pzc}{m}{it}%
             {<-> s * [1,150] pzcmi7t}{}
\DeclareMathAlphabet{\mathpzc}{OT1}{pzc}%
                                 {m}{it}
\input cyracc.def
\DeclareFontFamily{U}{russian}{}
\DeclareFontShape{U}{russian}{m}{n}
        { <5><6> wncyr5
        <7><8><9> wncyr7
        <10><10.95><12><14.4><17.28><20.74><24.88> wncyr10 }{}
\DeclareSymbolFont{Russian}{U}{russian}{m}{n}
\DeclareSymbolFontAlphabet{\mathcyr}{Russian}
\makeatletter
\let\@math@cyr\mathcyr
\renewcommand{\mathcyr}[1]{\@math@cyr{\cyracc #1}}
\makeatother

\renewcommand{\leq}{\leqslant}
\renewcommand{\geq}{\geqslant}
\renewcommand{\mathscr}{\mathcal}

\newcommand{\DM}{\mathit{DM}}

\newcommand{\DMe}{\smash{D\mspace{-2.mu}M^{\mathit{eff}}}}
\newcommand{\DMgme}{\smash{D\mspace{-2.mu}M^{\mathit{eff}}_{gm}}}

\newcommand{\Sp} [1]{\smash{\mathrm{Sp}^{\mathfrak S}_{#1}}}
\newcommand{\TDer} {\mathsf T}
\newcommand{\ste}{\mathcal E}

\renewcommand{\varinjlim}{\limind}%
\def\limind{\mathop{\oalign{\rm lim\cr
\hidewidth$\longrightarrow$\hidewidth\cr}}}%
\def\TO#1{\mathrel{\hbox to #1pt{\rightarrowfill}}}
\def\OT#1{\mathrel{\hbox to #1pt{\leftarrowfill}}}
\def\hocolim{\mathop{\oalign{\rm holim\cr
\hidewidth$\TO{26.5}$\hidewidth\cr}}}%
\newcommand{\Mod}{{\mathrm{\text{-}\mathrm Mod}}}
\newcommand{\derL}{\mathbf{L}}
\newcommand{\derR}{\mathbf{R}}
\newcommand{\sHom}{{\mathbf{Hom}}}

\newcommand{\ab} {{\mathit{Ab}}}

\newcommand{\smc}[1]{\mathscr Sm^{cor}_{#1}}

\newcommand{\ftr}[1] {\mathscr N^{tr}_{#1}}
\newcommand{\spt}[1] {Sp_{\mathrm{Tate}}(#1)}

\newcommand{\op}[1]{{ #1 }^{\mathit{op}}}

\newcommand{\Cyl}{{\mathrm{Cyl}}}

\newcommand{\sseq}[1]{{\mathscr{#1}^{\mathfrak S}}}

\newcommand{\Hom}{{\mathrm{Hom}}}

\newcommand{\Comp}{\mathit{Comp}}
\newcommand{\Htp} {\mathit K}
\newcommand{\Der} {\mathit D}
\newcommand{\ho} {\mathit{Ho}}

\newcommand{\unit}{\mathbf 1}

\newcommand{\To}{\longrightarrow}

\newcommand{\s}{\sigma}
\newcommand{\A}{\mathscr A}

\newcommand{\C}{\mathscr C}
\newcommand{\D}{\mathscr D}

\renewcommand{\L}{\mathscr L}

\newcommand{\R}{\mathscr R}

\newcommand{\T}{\mathscr T}

\newcommand{\X}{\mathscr X}

\newcommand{\G}{\mathcal G}
\renewcommand{\H}{\mathcal H}
\newcommand{\Ring}{\mathcal{R}}
\newcommand{\Site}{\mathcal{S}}

\newcommand{\NN} {\mathbb N}
\newcommand{\ZZ} {\mathbb Z}

\renewcommand{\AA} { \mathbb A }

\newcommand{\GG} { \mathbb G }

\begin{document}

\begin{abstract}
We define model category structures on the category of chain complexes
over a Grothendieck abelian category depending on the choice of a generating
family, and we study their behaviour with respect to tensor products
and stabilization. This gives convenient tools to construct and understand
triangulated categories of motives and we consider here the case of mixed motives 
over a regular base scheme.
\end{abstract}
\maketitle

\tableofcontents
\section*{Introduction}
This paper is an attempt to apply the general technics of abstract
homotopical algebra to describe derived categories of Grothendieck
categories. Our motivation
is to describe model categories of complexes of sheaves
``as locally as possible'', so that we can get easily
total left derived functors (like the tensor product).
The formalism we propose is made to stay close to
something that looks like descent theory.
Even though the results proved here are just applications
of a general theory (Bousfield localization, symmetric spectra),
our contribution consists to describe
these model structures with an emphasis on the
specific properties in the setting of Grothendieck categories.
As an example of properties that are not known for general
model categories, we give explicit descriptions of fibrations after Bousfield localizations
and prove that passing to symmetric spectra can preserve
the monoid axiom of Schwede and Shipley. We also prove that
any derived category of a Grothendieck category
can be obtained as a localization of the derived category of
modules over an additive category. This is used to give a sufficient condition
for the derived category of a Grothendieck category to be
compactly generated.

As an illustration, we give the construction of the categories
of motivic complexes and motivic spectra over a regular base scheme together
with their natural structure~: derived inverse and direct image functors, derived
tensor product, \emph{twisted exceptional direct image functors} for smooth morphisms.
In this construction, smooth always means abusively smooth of finite type.
This construction will be amplified and generalized 
in the forthcoming paper \cite{DMCD}. The results of this paper are also widely used
to understand the notion of mixed Weil cohomology,
which is introduced and studied in \cite{kunneth}.


\section{Model category structures on complexes}

\begin{paragr}\label{defdnsn}
Let $\A$ be an abelian category.
If $E$ is an object of $\A$ and $n$ an integer, we define $S^nE$
(resp. $D^nE$) as the acyclic complex concentrated in degree $n$ with $(S^nE)^{n}=E$
(resp. in degrees $n$ and $n+1$ with $(D^nE)^n=(D^nE)^{n+1}=E$
and as only non trivial differential the identity of $E$).
We then have canonical inclusions $S^{n+1}E\To D^nE$ induced by the identity of $E$.
The complex $D^{n}E$ represents the functor
$$\op{\Comp(\A)}\To\ab\quad ,\qquad C\longmapsto\Hom_\A(E,C^n) \ ,$$
and the complex $S^nE$ represents the functor
$$\op{\Comp(\A)}\To\ab\quad ,\qquad C\longmapsto\Hom_\A(E,Z^nC) \ ,$$
where $Z^nC=\ker(d:C^n\To C^{n+1})$ denotes the object of $n$-cocycles of the complex $C$.
We shall refer ourselves to \cite{Hov} for general
results and definitions on model categories, and to \cite{Hir}
for the definition of cellular model categories and for results
about Bousfield localizations.
\end{paragr}

\begin{thm}[Beke]\label{cmf1}
Let $\A$ be a Grothendieck abelian category.
The category $\Comp(\A)$ is a proper cellular model category with
quasi-isomorphisms as weak equivalences and the monomorphisms as
cofibrations.
\end{thm}
\begin{proof}
See \cite[Proposition 3.13]{Beke} for the construction
of this model category structure.
The fact that this model category is cellular comes from the fact that any
monomorphism in a Grothendieck category is effective, and that any
object in a Grothendieck category is compact (in fact accessible
according to \cite{SGA4}). Remark that fibrations are in particular
epimorphisms. This is proved as follows.
For any object $E$ in $\A$ and any integer $n$, the complex $D^nE$
is acyclic. Hence any fibration has the right lifting
property with respect to the inclusions $0\To D^{n}E$.
This lifting property for a map $p:C\To D$ is equivalent to say that for any
object $E$ of $\A$ and any integer $n$, the map
$\Hom_{\A}(E,C^{n})\To\Hom_{\A}(E,D^{n})$ is surjective.
This implies in particular that $p$ is degreewise split, hence
that $p$ is an epimorphism.
Properness now comes from classical homological algebra
as follows. Let
$$\xymatrix{
A\ar[r]^{i}\ar[d]_{p}&B\ar[d]^{q}\\
C\ar[r]_{j}&D}$$
be a pushout (resp. a pullback) square of complexes in $\A$.
Suppose that $i$ (resp. $q$) is a cofibration (resp. a fibration).
As $i$ (resp. $q$) is a monomorphism (resp. an epimorphism),
the same is true for $j$ (resp. for $p$). Moreover, it is obvious
that the cokernels of $i$ and $j$ (resp. the kernels of $p$ and $q$)
are canonically isomorphic, which implies we have a long exact
sequence of shape
$$H^{n}(A)\To H^{n}(B)\oplus H^{n}(C)\To H^{n}(D)\To H^{n+1}(A) \ .$$
We deduce then immediately that if $p$ (resp. $j$) is a quasi-isomorphism,
so is $q$ (resp. $i$), which is what we wanted to prove.
\end{proof}

\begin{num}
The above model category structure
will be called the \emph{injective model structure}.
The fibrations of this model structure will be called
the \emph{injective fibrations}, and the fibrant objects
the \emph{injectively fibrant objects}.\\
\indent The injective model structure gives a powerful tool to
construct right derived functors (well defined on
unbounded complexes). To define left derived functors
on unbounded complexes, we have to consider other model structures
depending on the choice of a set of generators of $\A$.
We shall restrict ourselves to a special situation
that will occur in the cases we are interested in.
This can be well formulated using the intuition of
descent theory.
\end{num}

\begin{df}
\label{descent_struct}
Let $\A$ be a Grothendieck category.\\
\indent If $\G$ is an essentially small set of objects of $\A$,
we define a morphism in $\Comp(\A)$ to be a \emph{$\G$-cofibration}
if it is contained in the smallest class of maps in $\Comp(\A)$
stable by pushouts, transfinite compositions and retracts
generated by the inclusions $S^{n+1}E\To D^nE$ for any integer $n$
and any $E$ in $\G$. A chain complex $C$ is \emph{$\G$-cofibrant}
if the map $0\To C$ is a $\G$-cofibration. For example,
If $C$ is a bounded above complex such that for any integer $n$,
$C^n$ is a direct sum of elements of $\G$, then $C$
is $\G$-cofibrant.\\
\indent A chain complex $C$ in $\A$
is \emph{$\G$-local}, or simply \emph{local} if no confusion
can arise from the context, if for any $E$ in $\G$
and integer $n$, the canonical map
$$\Hom_{\Htp(\A)}(E[n],C)\To\Hom_{\Der(\A)}(E[n],C)$$
is an isomorphism.\\
\indent If $\H$ is a small family of complexes in $\A$,
an object $C$ of $\Comp(\A)$ is \emph{$\H$-flasque}
(or simply \emph{flasque}) if for any integer $n$ and any $H$ in $\H$, the
abelian group $\Hom_{\Htp(\A)}(H,C[n])$ vanishes.\\
\indent A \emph{descent structure} on a Grothendieck category $\A$
is a couple $(\G,\H)$ where $\G$ is an essentially small
set of generators of $\A$, and $\H$ an essentially small set
of $\G$-cofibrant acylic complexes such that any $\H$-flasque
complex is $\G$-local.
\end{df}

\begin{ex}\label{exdesc}
Let $\Site$ be an essentially small Grothendieck site and $\Ring$ a sheaf of rings
on $\Site$. Then the category $\Ring\Mod$ of sheaves of $\Ring$-modules
is a Grothendieck category. For an object $X$ of $\Site$,
denote by $\Ring(X)$ the free sheaf of $\Ring$-modules generated by $X$.
The family $\G_\Site$ of the $\Ring(X)$'s form a generating set for $\Ring\Mod$.
If $\X$ is an hypercover of $X$, then the simplicial $\Ring$-module
$\Ring(\X)$ can be seen as a complex of $\Ring$-modules with the differentials
given by the alternated sums of the face operators.
For an hypercover $\X$ of an object $X$, denote by $\widetilde{\Ring}(\X)$
the cone of the structural map $\Ring(\X)\To\Ring(X)$.
The complexes $\Ring(X)$, $\Ring(\X)$ and $\widetilde{\Ring}(\X)$ are
then $\G_\Site$-cofibrant.
Let $\H_\Site$ be the family of all the complexes of the form
$\widetilde{\Ring}(\X)$ for any hypercover $\X$ of an object $X$ of $\Site$.
Then the couple $(\G_\Site,\H_\Site)$ is a descent structure on the
category $\Ring\Mod$. One way to prove this comes from
Verdier's computation of hypercohomology; see \cite[Expos\'e V]{SGA4}
and \cite{Br}. If $X$ is an object of $\Site$ and $C$ is a complex of $\R$-modules,
then we have the formula
$$\limind_{\X\in H_{X}}\Hom_{\Htp(\Ring\Mod)}(\R(\X),C)=\Hom_{\Der(\Ring\Mod)}(\R(X),C)$$
where $H_{X}$ denotes the filtering category of hypercoverings of $X$ up to simplicial
homotopy (and $\Hom_{\Der(\Ring\Mod)}(\R(X),C)$ is nothing else but the hypercohomlogy
of $X$ with coefficients in $C$). We deduce that if $C$ is $\H_\Site$-flasque, 
as we have distinguished triangles in $\Htp(\Ring\Mod)$ of shape
$$\R(\X)\To\R(X)\To\widetilde{\Ring}(\X)\To\Ring(\X)[1] $$
for any $\X$ in $H_{X}$, we have isomorphisms
$$\Hom_{\Htp(\Ring\Mod)}(\R(X),C)=\Hom_{\Htp(\Ring\Mod)}(\R(\X),C) \  .$$
As $H_{X}$ if filtering, we deduce that the canonical map
$$\Hom_{\Htp(\Ring\Mod)}(\R(X),C)\To\limind_{\X\in H_{X}}\Hom_{\Htp(\Ring\Mod)}(\R(\X),C)$$
is an isomorphism. Hence $C$ is $\G_\Site$-local.
\end{ex}

\begin{ex}\label{exshtrans}
Let $S$ be a regular scheme,
 and $\ftr S$ be the category of Nisnevich sheaves with transfers over $S$
defined in \cite[4.2.3]{Deg7}.
Let $\G_S$ be the category of sheaves with transfers 
of type $L_S(X)$ for any smooth $S$-scheme $X$.
Then it follows from proposition 4.2.8 of \emph{loc. cit.} that
$\A_S$ is a Grothendieck abelian category and that $\G_S$
is a set of generators. We denote by
$\H_S$ the set of complexes obtained as the cones
of the maps $L_S(\X)\To L_S(X)$, where
$X$ is any smooth $S$-scheme
and $\X$ is any Nisnevich hypercovering of $X$.
It follows from the example above and from
propositions 4.2.5 and 4.2.9 of \emph{loc. cit.}
that $(\G_S,\H_S)$ is a descent structure on $\ftr S$.
\end{ex}

\begin{thm}\label{cmf2}
Let $\A$ be a Grothendieck category endowed with a
descent structure $(\G,\H)$.
Then the category $\Comp(\A)$ is a proper cellular model category
with quasi-isomorphisms as weak equivalences,
and the $\G$-cofibrations as cofibrations.
Furthermore, for a complex $C$ in $\A$, the following conditions
are equivalent.
\begin{itemize}
\item[(i)] $C$ is fibrant;
\item[(ii)] $C$ is $\H$-flasque;
\item[(iii)] $C$ is $\G$-local.
\end{itemize}
This model structure on the category $\Comp(\A)$
will be called the \emph{$\G$-model structure}.
\end{thm}

\begin{rem}\label{remcarfibgloc}
It is possible to give a nice description of the
fibrations for this model category structure; see Corollary \ref{characfibrations2}.
\end{rem}

\begin{num}\label{defJcof}
The proof of this theorem will be completed
by the results below. We first have to
fix some notations.\\
\indent We define $I$ to be the set of inclusions of the form
$S^{n+1}E\To D^nE$ for any integer $n$ and any $E$ in $\G$.\\
\indent If $C$ is a complex, $\Cyl(C)$ is the complex
defined by
$$\Cyl(C)^n=C^n\oplus C^{n+1}\oplus C^n$$
with differential $d(x,y,z)=(dx-y,-dy,y+dz)$.
The projection $\s :\Cyl(C)\To C$, defined by the formula
$\s(x,y,z)=x+z$, is a quasi-isomorphism.
We also have a canonical inclusion $(i_0,i_1):C\oplus C\To \Cyl(C)$
defined by $i_0(x)=(x,0,0)$ and $i_1(z)=(0,0,z)$, so that
the composite $\s(i_0,i_1)$ is equal to the codiagonal
from $C\oplus C$ to $C$. It is an easy exercise to check that to define
a chain homotopy between two maps of complexes
$u_0$ and $u_1$ from $C$ to $C'$ is equivalent to define a map $h$
from $\Cyl(C)$ to $C'$ such that $hi_e=u_e$ for $e=0,1$.\\
\indent We define at last $J=J'\cup J''$, where $J'$ is the set of maps
$0\To D^nE$ for any integer $n$, any $E$ in $\G$,
and $J''$ is the set of maps $H\oplus H[n]\To \Cyl(H)[n]$
for any $H$ in $\H$ and any integer $n$.\\
\indent To prove theorem \ref{cmf2}, we are going to check the conditions
of \cite[Theorem 2.1.19]{Hov}. We adopt the notations of
\cite{Hov}: $I$-$\mathit{cof}$ (resp. $J$-$\mathit{cof}$) is the class of $\G$-cofibrations
(resp. the smallest class of maps stable by pushouts, transfinite
compositions and retracts that contains $J$), and $I$-$\mathit{inj}$ (resp. $J$-$\mathit{inj}$)
the class of maps with the right lifting property with respect
to $I$ (resp. $J$). We will prove that $I$-$\mathit{cof}$, $J$-$\mathit{cof}$,
$I$-$\mathit{inj}$ and $J$-$\mathit{inj}$ are respectively the classe
of cofibration, the class of trivial cofibrations, the class of trivial fibrations
and the class of fibrations for this model structure.
\end{num}

\begin{lm}\label{Icof}
For any $E$ in $\G$ and any integer $n$, the complexes
$S^{n}E$ and $D^{n}E$ are $\G$-cofibrant.
For any $\G$-cofibrant complex $C$, the inclusion
$C\oplus C\To \Cyl(C)$ is a $\G$-cofibration. If
furthermore $C$ is acyclic, then this map is also
a quasi-isomorphism.
\end{lm}

\begin{proof}
The pushout square
$$\xymatrix{
S^{n}E\ar[r]\ar[d]&0\ar[d]\\
D^{n-1}E\ar[r]&S^{n+1}E}$$
implies the first assertion.

Remember that the cone $\mathrm{Cone}(p)$ of a morphism 
of complexes $p:X\To Y$ is defined by
$\mathrm{Cone}(p)^n=Y^n\oplus X^{n+1}$, with
differentials given by the formula $d(y,x)=(dy+px,-dx)$.
We then have an obvious short exact sequence
$$
\xymatrix{0\ar[r]&Y\ar^/-6pt/u[r]&{\mathrm{Cone}}(p)\ar^/3pt/v[r]&X[1]\ar[r]&0}
$$
where $u(y)=(y,0)$, and $v(y,x)=x$.

Consider an element $E$ of $\G$ and an integer $n$.
One checks $\mathrm{Cone}(1_{S^{n+1}E})=D^nE$ and 
the inclusion of $S^{n+1}E$ in $\mathrm{Cone}(1_{S^{n+1}E})$
is the canonical inclusion of $S^{n+1}E$ in $D^nE$.
The complex $\mathrm{Cone}(1_{D^nE})$ is the split acyclic complex
$$0\To E\To E\oplus E\To E\To 0$$
where $E\oplus E$ is placed in degree $n+1$.
One deduces that the canonical map
$$D^nE\amalg_{S^{n+1}E}\mathrm{Cone}(1_{S^{n+1}E})\To\mathrm{Cone}(1_{D^nE})$$
is a $\G$-cofibration: it corresponds to the obvious inclusion of the complex
$$0\To 0\To E\oplus E\To E\To 0$$
in $\mathrm{Cone}(1_{D^nE})$.
In particular, the map $\mathrm{Cone}(1_{S^{n+1}E})\To\mathrm{Cone}(1_{D^nE})$
is a $\G$-cofibration.
Let $\mathbf{C}$ be the class of maps $A\To B$ such that
the map
$$B\amalg_A\mathrm{Cone}(1_A)\To\mathrm{Cone}(1_B)$$
induced by the commutative square
$$\xymatrix{A\ar[r]\ar[d]&{\mathrm{Cone}}(1_A)\ar[d]\\
B\ar[r]&{\mathrm{Cone}}(1_B)}$$
is a $\G$-cofibration. This class of maps is stable
by pushouts, transinite compositions and retracts (exercise left to the reader).
As $\mathbf{C}$, contains $I$, it has to contain all of $I$-$\mathit{cof}$.
Hence if $C$ is $\G$-cofibrant, as the map $0\To C$ is in $\mathbf{C}$,
the map $C\To\mathrm{Cone}(1_C)$ is a $\G$-cofibration.
We deduce that
for any $\G$-cofibrant complex $C$, the canonical inclusion
$C\To\mathrm{Cone}(1_C)$ is a $\G$-cofibration.
Indeed, we have a pushout square of the following form
$$\xymatrix{
&C\ar[r]\ar[d]&C\oplus C\ar[d]&\\
&{\mathrm{Cone}}(1_C)\ar[r]&\Cyl(C)&}$$
and the $\G$-cofibrations are stable by pushout.

If moreover $C$ is acyclic, then so is $C\oplus C$.
As the map $\Cyl(C)\To C$ is always a quasi-isomorphism,
the complex $\Cyl(C)$ is acyclic as well. Hence
the inclusion $C\oplus C\To\Cyl(C)$ has to be a
quasi-isomorphism.
\end{proof}

\begin{lm}\label{JdansIW}
Any element of $J$-$\mathit{cof}$ is both a $\G$-cofibration
and a quasi-iso\-mor\-phism.
\end{lm}

\begin{proof}
Any $\G$-cofibration is a monomorphism
and the class of maps
which are mono\-mor\-phisms and quasi-isomorphisms
is stable by pushouts, transfinite compositions and retracts.
Hence the class of $\G$-cofibrations which are
quasi-isomorphisms is also stable by any of these
operations. It is then sufficient to check that any map in $J$ is
a $\G$-cofibration and a quasi-isomorphism.
This follows from lemma \ref{Icof}.
\end{proof}

\begin{lm}\label{Jfibfl}
Let $C$ be a complex in $\A$.
Then the map $C\To 0$ is in
$J$-$\mathit{inj}$ if and only if $C$ is $\H$-flasque.
\end{lm}

\begin{proof}
Let $H$ be in $\H$, and $u$ a map
from $H$ to $C[n]$. The data of a chain
homotopy between $u$ and the zero map
is equivalent to the data of a lifting $k$
in the commutative square below.
$$\xymatrix{
H\oplus H\ar[r]^{(u,0)}\ar[d]&C[n]\ar[d]\\
\Cyl(H)\ar[r]\ar@{..>}[ur]^k&0}$$
As the right lifting property of $C\To 0$
with respect to the maps of the shape $0\To A$
for any complex $A$ is always verified, this proves
the assertion.
\end{proof}

\begin{lm}\label{Hflacycl}
An $\H$-flasque complex $C$ is acyclic if and only if
the map $C\To 0$ has the right lifting property with respect to $I$.
\end{lm}

\begin{proof}
For any complex $C$ in $\A$,
as the cohomology group $H^n(\Hom_\A(E,C))$ is the cokernel
of the map
$$\Hom_{\Comp(\A)}(D^{n-1}E,C)\To\Hom_{\Comp(\A)}(S^{n}E,C) \ , $$
it is clear that $C\To 0$ is in $I$-$\mathit{inj}$ if and only if
the complexes $\Hom_\A(E,C)$ are acyclic for any $E$ in $\G$.
If $C$ is $\H$-flasque, then for any $E$ in $\G$ and integer $n$,
we have canonical isomorphisms
$$H^n(\Hom_\A(E,C))=\Hom_{\Htp(\A)}(E,C[n])=\Hom_{\Der(\A)}(E,C[n]) \ .$$
To finish the proof, it is then sufficient to say that a complex $C$
in $\A$ is acyclic if and only if for any $E$ in $\G$ and integer $n$,
the group $\Hom_{\Der(\A)}(E,C[n])$ vanishes.
\end{proof}

\begin{lm}\label{carjinjqis}
A map in $J$-$\mathit{inj}$ is a quasi-isomorphism if and only if
it is in $I$-$\mathit{inj}$.
\end{lm}

\begin{proof}
Let $p:X\To Y$ be map in $J$-$\mathit{inj}$. We can see that it is an epimorphism:
for any element $E$ of $\G$ and any integer $n$, the map
$$\Hom_\A(E,X^n)=\Hom_{\Comp(\A)}(D^nE,X)\to\Hom_{\Comp(\A)}(D^nE,Y)=\Hom_\A(E,Y^n)$$
is surjective. As $\G$ is a set of generators of $\A$, this implies that $p$
is an epimorphism. In particular $\ker p$ is acyclic if and only if
$p$ is a quasi-isomorphism. If $p$ is in $I$-$\mathit{inj}$,
the map $\ker p\To 0$ is also in $I$-$\mathit{inj}$.
But $J\subset$~$I$-$\mathit{cof}$ (see \ref{JdansIW}),
which implies by lemma \ref{Jfibfl}
that $\ker p$ is $\H$-flasque. Hence by lemma \ref{Hflacycl},
the complex $\ker p$ is acyclic, and $p$ is a quasi-isomorphism.

Consider now a quasi-isomorphism $p$ in $J$-$\mathit{inj}$.
As $p$ is then an epimorphism, $\ker p$ is acyclic,
and as right lifting properties are stable by pull-back, the map $\ker p\To 0$
is in $J$-$\mathit{inj}$. We deduce from \ref{Jfibfl}
and \ref{Hflacycl} that the map $\ker p\To 0$ is in
$I$-$\mathit{inj}$. Let
$$\xymatrix{S^{n+1}E\ar[r]^x\ar[d]&X\ar[d]^p\\D^nE\ar[r]_y&Y}$$
be a commutative square in $\Comp(\A)$ where $E$ is in $\G$, and $n$
is an integer. It corresponds to maps $x:E\To X^{n+1}$ and
$y:E\To Y^{n}$ such that $px=dy$.
To show that the square above admits a lift is equivalent to
show that there is a map $\xi:E\To X^n$ such that $d\xi=x$ and $p\xi=y$.
As $p$ has the right lifting property with respect to $J$,
$\Hom_\A(E,p^n)$ is surjective, and so there exists a map
$x':E\To X^n$ such that $px'=y$. We know that $pdx'=dpx'=dy=px$,
so that $x-dx'$ factors through $\ker p$. As the map $\ker p\To 0$
has the right lifting property with respect to $I$, this implies that
there is a map $x'':E\To X^n$ such that $px''=0$ and $dx''=x-dx'$.
We can know set $\xi=x'+x''$. We have that $p\xi=px'+px''=y$
and $d\xi=dx'+dx''=x$. Hence the desired lift.
\end{proof}

\begin{proof}[Proof of theorem \ref{cmf2}]
One deduces from \ref{JdansIW} and \ref{carjinjqis}
that all the conditions of \cite[Theorem 2.1.19]{Hov}
are satisfied. In particular, we obtain that
we have defined a cofibrantly generated
model category with quasi-isomorphisms
as weak equivalences and with $I$ (resp. $J$)
as a set of generators for the cofibrations
(resp. for the trivial cofibrations). To see
that this structure is cellular, it is sufficient to say that
any $\G$-cofibration is a monomorphism
and that any monomorphism in a Grothendieck
category is effective. Properness comes from the fact
that cofibrations (resp. fibrations) are monomorphisms
(resp. epimorphisms) and standard results of
homological algebra (see the end of the proof of \ref{cmf1}).\\
\indent It remains to prove the equivalence of conditions
(i)-(iii). Using definitions and \ref{Jfibfl}, it is sufficient
to show that (iii) implies (ii). Let $C$ be a $\G$-local
complex. One can prove that the full
subcategory of $\Htp(\A)$ made of $\G$-cofibrant
complexes is contained in the localizing subcategory
of $\Htp(\A)$ generated by $\G$ (in the sense of \cite{Nee}).
Then for any $\G$-cofibrant complex $A$,
the map below is an isomorphism.
$$\Hom_{\Htp(\A)}(A,C)\To\Hom_{\Der(\A)}(A,C)$$
In particular, for any $H$ in $\H$ and integer $n$, we have
$\Hom_{\Htp(\A)}(H,C[n])=0$. Hence $C$ is $\H$-flasque.
\end{proof}

\begin{rem}\label{remabstractnonsense}
The model structure of theorem \ref{cmf2}
depends only on $\G$ and not on the whole
descent structure $(\G,\H)$. The data $\H$
is only a tool to get a nice characterization of the
fibrant objects in terms of flasque properties (as we will
see below, that will be usefull to construct total left derived functors).
One can prove directly by abstract non sense that for
any Grothendieck category $\A$ with a choice of
an essentially small set of generators $\G$,
there is Quillen model category structure on $\Comp(\A)$
with quasi-isomorphisms as weak equivalences
and $\G$-cofibrations as cofibrations.
To show this, just remark that the class of
quasi-isomorphisms is accessible (see the proof of
\cite[proposition 3.13]{Beke}), that any map with
the right lifting property with respect to $\G$-cofibrations
is a quasi-isomorphism (by virtue of \cite[Proposition 1.4]{Hov2}),
and apply Jeff Smith's Theorem (see \cite[Theorem 1.7 and Proposition 1.15]{Beke}).
In fact one can show that any small set $\G$
of generators of $\A$ can be completed in a
descent structure for $\A$, in such a way that the description
of this model category structure given by \ref{cmf2} is always available
(but in practice, we essentially always have a canonical choice of
descent structure). One can proceed as follows.
Chose a set $J_{0}$ of generating trivial
cofibration (with respect to the model structure above).
Define $\H$ to be the set of the cokernels of
the maps in $J_{0}$. Then $(\G,\H)$ is a descent structure.
\end{rem}

\begin{ex}\label{Dftr}
Consider the notations of \ref{exshtrans}.

Let $C$ be a complex of sheaves with transfers over $S$.

For any smooth $S$-scheme $X$, we denote by $H^n C(X)$
 (resp. $H^n(X,C)$)
the $n$-th cohomology of the complex of abelian groups 
made by the sections of $C$ over $X$ (resp. the $n$-th Nisnevich
hypercohomology of $X$ with coefficients in the 
complex $C$ where we forget transfers).

We will say that $C$ is Nisnevich fibrant if it is fibrant
with respect to the descent structure $(\G_S,\H_S)$.
According to the previous theorem and \cite[4.2.9]{Deg7},
this means one of the following equivalent conditions is true~:
\begin{itemize}
\item[(i)] For any smooth $S$-scheme $X$ and any integer $n$,
 the canonical map $H^n(C(X)) \rightarrow H^n_{\mathrm{Nis}}(X;C)$
  is an isomorphism.
\item[(ii)] For any smooth $S$-scheme $X$,
 any Nisnevich hypercovering $\X$ of $X$ and any integer $n$, 
  the canonical map $H^n(C(X)) \rightarrow H^n(\mathrm{Tot}\,C(\X))$ is an isomorphism.
\end{itemize}
\end{ex}

\begin{num}\label{defstdesc}
Let $\A$ and $\A'$ be two Grothendiek categories.
Suppose that $(\G,\H)$ (resp. $(\G',\H')$) is
a descent structure on $\A$ (resp. $\A'$).
A functor $f^*:\A'\To\A$ \emph{satisfies descent} (with respect to the
above descent structures) if it satisfies the following conditions.
\begin{itemize}
\item[(i)] The functor $f^*$ preserves colimits.
\item[(ii)] For any $E'$ in $\G'$,
$f^*(E')$ is a direct sum of elements of $\G$.
\item[(iii)] For any $H'$ in $\H'$, $f^*(H')$
is in $\H$.
\end{itemize}
It is a standard fact that a functor
between Grothendieck categories preserves
colimits if and only if it has a right adjoint.
\end{num}

\begin{thm}\label{fonct1}
Under the assumptions above,
let $f^*:\A'\To\A$ be a functor that preserves colimits,
and $f_*:\A\To\A'$ be its right adjoint.
If $f^*$ satisfies descent, then the pair of adjoint functor
$$f^*:\Comp(\A')\To\Comp(\A)\quad\text{and}
\quad f_*:\Comp(\A)\To\Comp(\A')$$
is a Quillen pair with respect to the model
structures associated to $\G$ and $\G'$
by the theorem \ref{cmf2}. In particular,
the functors $f^*$ and $f_*$ have the functors
$$\derL f^*:\Der(\A')\To\Der(\A)\quad\text{and}\quad
\derR f_*:\Der(\A)\To\Der(\A')$$
as left and right derived functors
respectively, and $\derL f^*$ is left adjoint
to $\derR f_*$.
\end{thm}

\begin{proof}
If $f^*$ satisfies descent, then one checks
easily that it preserves cofibrations; see \cite[Lemma 2.1.20]{Hov}.
It is then sufficient to show that it sends
trivial cofibrations to quasi-isomorphisms.
As $f^*$ preserves colimits and cofibrations,
by virtue of \emph{loc. cit}, it is sufficient to prove that any element of $f^*(J')$
is a quasi-isomorphism, where $J'$
is the set of trivial cofibrations defined as in \ref{defJcof}.
This comes directly from the fact that
$f^*(\H')$ is contained in $\H$
and that $f^*$ commutes to the formation
of cones and cylinders (just
because it is additive).
\end{proof}

\begin{rem}\label{composefctder}
Consider three Grothendieck categories
endowed with descent structures $\A$, $\A'$, and $\A''$.
Let $f^{\prime \, *}:\A''\To\A'$ and $f^*:\A'\To\A$
be two functors commuting to colimits
and satisfying descent, and denote by $f'_*$
and $f_*$ their corresponding right adjoints.
Then it follows easily from
general abstract nonsense about Quillen
adjunctions and the preceding theorem
that we have canonical isomorphisms
of total derived functors
$$\derL f^{*}\circ\derL f^{\prime \, *}\simeq\derL(f^{*}\circ f^{\prime\, *})\quad\text{and}
\quad \derR(f'_*\circ f_*)\simeq\derR f'_*\circ\derR f_* \ .$$
\end{rem}

\begin{ex} \label{DM_4functors}
Consider the notations of Example \ref{exshtrans}.

Let $f:T \rightarrow S$ be a morphism between regular schemes.
In \cite[4.2.5]{Deg7} is defined the base change functor
$f^*:\ftr S \rightarrow \ftr T$ together with its right adjoint
$f_*:\ftr T \rightarrow \ftr S$. By definition, 
 $f^* L_S(X)=L_T(X \times_S T)$
  for any smooth finite type $S$-scheme $X$.
Thus $f^*$ obviously satisfies descent with respect to the structure
$(\G,\H)$ defined in \ref{exshtrans}. 
Thus we get a pair of adjoint derived functors
$$\derL f^*:\Der(\ftr S)\To\Der(\ftr T)\quad\text{and}\quad
\derR f_*:\Der(\ftr T)\To\Der(\ftr S).$$

Suppose moreover that $f:T \rightarrow S$ is smooth of finite 
type. Then in \cite[4.2.5]{Deg7} is defined the functor
$f_\sharp:\ftr T \rightarrow \ftr S$ which "forgets the base"~:
for any smooth finite type $T$-scheme $Y \xrightarrow \pi T$, 
$f_\sharp L_T(Y)=L_S(Y \xrightarrow \pi T \xrightarrow f S)$.
It is proven in \emph{loc. cit.} that this functor as for 
right adjoint the functor $f^*$. Obviously, $f_\sharp$ satisfies
descent with respect to the structure $(\G,\H)$, thus we get
a pair of adjoint derived functors
$$\derL f_\sharp:\Der(\ftr S)\To\Der(\ftr T)\quad\text{and}\quad
\derR f^*=f^*:\Der(\ftr T)\To\Der(\ftr S).$$
\end{ex}

\section{Derived tensor product}

\begin{num}\label{defflatfamily}
Let $\A$ be a Grothendieck category
endowed with a closed symmetric monoidal structure\footnote{Remember that
a closed symmetric monoidal category is a symmetric monoidal category with
internal $\Hom$'s.}. An object $X$ of $\A$ will be said \emph{flat}
if the functor defined by the tensor product with $X$
is exact. An essentially small family of generators $\G$ on $\A$
will be said to be \emph{flat} if the following conditions
are satisfied.
\begin{itemize}
\item[(i)] Any element of $\G$ is flat.
\item[(ii)] The set $\G$ is stable by tensor product up to
isomorphism and the unit object of $\A$ is in $\G$.
\end{itemize}
The category $\Comp(\A)$ of complexes in $\A$ is
canonically endowed with a symmetric monoidal structure
induced by the one of $\A$ by setting for two complexes $X$ and $Y$
$$(X\otimes Y)^n=\bigoplus_{p+q=n}X^p\otimes Y^q \ .$$
The differential of $X\otimes Y$ is given by the
Leibniz formula
$$d(x\otimes y)=dx\otimes y+(-1)^p x\otimes dy$$
where $p$ is the degree of $x$. The unit is just the unit of $\A$
concentrated in degree zero. The associativity structure
comes from the one on $\A$, and the symmetry is induced
by the one on $\A$ with the usual sign convention:
$$X\otimes Y\To Y\otimes X\quad , \qquad x\otimes y\longmapsto (-1)^{pq}y\otimes x$$
where $p$ and $q$ are the degree of $x$ and $y$ respectively.
A descent structure $(\G,\H)$ on $\A$ is \emph{weakly flat}
if the following conditions are satisfied.
\begin{itemize}
\item[(i)] The set $\G$ is stable by tensor product up to
isomorphism and the unit object of $\A$ is in $\G$.
\item[(ii)] For any $E$ in $\G$ and any $H$ in $\H$, the complex
$E\otimes H$ is acyclic.
\end{itemize}
A descent structure $(\G,\H)$ on $\A$ is \emph{flat}
if it is weakly flat, and if for any complex $C$ on $\A$ and any
$H$ in $\H$, the complex $C\otimes H$ is acyclic.
We will see that a descent structure is flat if and only if
the underlying generating family is flat; see Proposition \ref{cmf333} below.
\end{num}

\begin{ex}\label{exdescflat}
Consider a Grothendieck site $\Site$
and a sheaf of rings $\Ring$ on $\Site$.
Consider the usual tensor product
on the category of sheaves of $\Ring$-modules.
If the category $\Site$ has finite products, then
the descent structure of Example \ref{exdesc}
is flat. This comes from the fact that for two objects
$X$ and $Y$ in $\Site$, we have the formula
$\Ring(X\times Y)=\Ring(X)\otimes^{}_{\Ring}\Ring(Y)$
and from the fact that the free sheaves of shape $\Ring(X)$
are flat.
\end{ex}

\begin{prop}\label{carflat1}
Let $(\G,\H)$ be a weakly flat descent structure on $\A$.
The corresponding $\G$-model structure on $\Comp(\A)$
is a symmetric monoidal model category (see \cite{Hov}).
In particular, the tensor product on $\Comp(\A)$ has a total left
derived functor
$$\Der(\A)\times\Der(\A)\To\Der(\A)\quad ,\qquad (A,B)\longmapsto A\otimes^\derL B \ .$$
\end{prop}

\begin{proof}
First of all, consider two $\G$-cofibrations
$a:A\To A'$ and $b:B\To B'$. Then the map
$$c:A\otimes B'\amalg_{A\otimes B}A'\otimes B\To A'\otimes B'$$
induced by the commutative square
$$\xymatrix{A\otimes B\ar[r]\ar[d]&A\otimes B'\ar[d]\\
A'\otimes B\ar[r]&A'\otimes B'}$$
is a $\G$-cofibration: according to \cite[Lemma 4.2.4]{Hov},
it is sufficient to check this property when $a$ and $b$ are in $I$,
but in this case, this follows easily from condition (i) of the definition
of weak flatness.
If furthermore $b$ is a quasi-isomorphism, then
so is $c$: by \emph{loc. cit.}, it is sufficient to prove this
assuming that $a$ is in $I$ and $b$ is in $J$, in which case this is trivial
using condition (ii). This achieves the proof.
\end{proof}

\begin{ex}\label{DMmonoidal}
Consider the notations of Example \ref{exshtrans}.

The category $\ftr S$ admits a closed symmetric monoidal structure
from \cite[4.2.12, 4.12.14]{Deg7}. This structure satisfies the
fundamental property that $L_S(X) \otimes L_S(Y)=L_S(X \times_S Y)$.
Thus, the descent structure $(\G,\H)$ is obviously weakly flat and
we get a left derived tensor product
 $\otimes^\derL$ on $\Der(\ftr S)$.

However, we do not know if $(\G_S,\H_S)$ is flat - which actually leads 
us to consider the weakly flat condition.
\end{ex}

\begin{prop}\label{monoidaxiom}
Let $(\G,\H)$ be a weakly flat descent structure on a closed symmetric monoidal
Grothendieck category, and $C$ be a complex in $\A$.
We suppose that for any $H$ in $\H$, and for any integer $n$,
the complex $H\otimes C^n$ is acyclic.

Then the functor $A\longmapsto A\otimes C$
is a left Quillen functor from the $\G$-model structure on $\Comp(\A)$
to the injective model structure on $\Comp(\A)$. In particular, for any
complex $C$ over $\A$ and any quasi-isomorphism which is also a $\G$-cofibration $X\To Y$,
the map $C\otimes X\To C\otimes Y$ is a quasi-isomorphism and a monomorphism.
\end{prop}

\begin{proof}
For any $H$ in $\H$, the complex $H\otimes C$ is acyclic:
this a total complex of double complex whose rows $H\otimes C^n$
are exact.

To prove the tensor product by $C$ defines a left Quillen
functor, applying \cite[Lemma 2.1.20]{Hov},
it is sufficient to prove that for any map $f$ in $I$ (resp. in $J$),
the map $f\otimes C$ is a monomorphism (resp. a quasi-isomorphism).
The only non trivial part to check is that for any $H$ in $\H$,
the map
$$H\otimes C\oplus H\otimes C\To \Cyl(H)\otimes C$$
is a quasi-isomorphism, or in other words, that $\Cyl(H)\otimes C$
is acyclic. But we have a push-out square
$$\xymatrix{
H\otimes C\ar[r]\ar[d]&H\otimes C\oplus H\otimes C\ar[d]\\
{\mathrm{Cone}}(1_H)\otimes C\ar[r]&\Cyl(H)\otimes C}$$
where the left vertical map is a monomorphism.
As we already know that $H\otimes C$ is acyclic, it is sufficient
to check that $\mathrm{Cone}(1_H)\otimes C$ is acyclic.
If $C$ is concentrated in one degree, one has
$\mathrm{Cone}(1_H)\otimes C=\mathrm{Cone}(1_{H\otimes C})$,
and this is trivial. The general case follows as this implies
that $\mathrm{Cone}(1_H)\otimes C$
is the total complex associated to a double complex with acyclic rows.
\end{proof}

\begin{cor}\label{cmf3}
Let $\A$ be a closed symmetric monoidal Grothendieck category
endowed with a flat descent structure $(\G,\H)$.
Then the associated $\G$-model structure is a monoidal
model category satisfying the monoid axiom\footnote{The monoid
axiom is a technical property that makes life easier
when one wants to define the homotopy theory of monoids or of
modules over a given monoid in a symmetric monoidal model category.
In our particular case, this can be seen as the second statement
in proposition \ref{monoidaxiom}.} (see \cite{SS}).
\end{cor}

\begin{proof}
The first assertion follows from proposition \ref{carflat1}.
As the trivial cofibrations of the injective model
structure on $\Comp(\A)$ are in particular stable by
pushouts and transfinite compositions, the monoid axiom
follows from proposition \ref{monoidaxiom}.
\end{proof}

\begin{lm}\label{cmf33}
Let $\A$ be a closed symmetric monoidal Grothendieck category,
and $\G$ be an essentially small set of generators of $\A$.
Consider a map of complexes $i:A\To B$ such that for any
$E$ in $\G$, $i\otimes E$ is a quasi-isomorphism. Then for any
$\G$-cofibrant complex $C$, the map of complexes $i\otimes C$
is a quasi-isomorphism.
\end{lm}

\begin{proof}
One checks easily that the map of complexes $i\otimes S^nE$ and $i\otimes D^nE$
are quasi-isomorphisms for any $E$ in $\G$ and any integer $n$.
If $C'$ is a direct summand of $C$, and if $i\otimes C$
is a quasi-isomorphism, then so is $i\otimes C'$. If
$$\xymatrix{S^{n+1}E\ar[r]\ar[d]&C\ar[d]\\D^nE\ar[r]&C'}$$
is a pushout in $\Comp(\A)$ with $E$ in $\G$ and $n$ an integer,
then if $i\otimes C$ is acylic, so is $i\otimes C'$: this comes from
the cube lemma applied to the injective model structure (see \cite[Lemma 5.2.6]{Hov}),
and the facts that the induced map
$A\otimes S^nE\To A\otimes D^nE$ is a monomorphism for any $A$,
and that the tensor product by any complex preserves pushouts (remember it
has a right adjoint). We know by the small object argument applied to $I$
that any $\G$-cofibrant complex $C$ is a direct summand
of a well ordered colimit of complexes constructed inductively
by some pushout of the type above. As quasi-isomorphisms are stable by
filtered colimits in a Grothendieck category, one deduces
that for any $\G$-cofibrant complex $C$, $i\otimes C$ is a quasi-isomorphism.
\end{proof}

\begin{prop}\label{cmf333}
Let $\A$ be a closed symmetric monoidal Grothendieck category,
and $(\G,\H)$ a descent structure on $\A$.
If $\G$ is flat, then the descent structure $(\G,\H)$ is flat.
Moreover, the tensor product by any $\G$-cofibrant
complex preserves quasi-isomorphisms.
\end{prop}

\begin{proof}
Suppose that $\G$ is flat. It is then clear that the descent structure is weakly flat.
In particular, the tensor product on $\Comp(\A)$ has a left derived functor
defined using the $\G$-model structure (\ref{carflat1}).
Let $C$ be a $\G$-cofibrant complex. We can see by lemma \ref{cmf33}
that the functor $X\longmapsto C\otimes X$ preserves quasi-isomorphisms,
so that for any complex $X$ over $\A$, the canonical map
$C\otimes^\derL X\to C\otimes X$
is an isomorphism in $\Der(\A)$. One concludes that for any
quasi-isomorphism between $\G$-cofibrant complexes $C\To C'$ and any
complex $X$, the map $C\otimes X\To C'\otimes X$ is a quasi-isomorphism.
In particular, for any $H$ in $\H$ and any complex $X$, $H\otimes X$
is acyclic and the descent structure $(\G,\H)$ is flat.
\end{proof}

\begin{rem}\label{tensdercoincide}
If $\G$ is a flat generating family of $\A$, then for any complex $X$,
the total left derived functor of $Y\longmapsto X\otimes Y$ (cf. \ref{monoidaxiom})
is isomorphic to the functor $Y\longmapsto X\otimes^{\derL} Y$
(obtained as the total left derived functor of the
functor $(A,B)\longmapsto A\otimes B$; cf. \ref{carflat1}).
This is an immediate consequence of the fact that
the tensor product by a $\G$-cofibrant
complex preserves quasi-isomorphisms.
This property is of course interesting by itself
but will be used in a very essential way to prove
the monoid axiom after further localizations (see e.g. Corolllary \ref{cmf4}).
\end{rem}

\begin{num}\label{defweaklyflatobject}
Consider a closed symmetric monoidal Grothendieck category $\A$
endowed with a weakly flat descent structure $(\G,\H)$.

Let $T$ be a given object of $\A$,
define $\G[T]$ to be the family of objects of the form $E\otimes T^{\otimes n}$
with $E$ in $\G$ and $n\geq 0$.

We will say that $T$ is \emph{weakly flat} (with respect to $(\G,\H)$)
if it has the two following properties:
\begin{itemize}
\item[(a)] for any $H$ in $\H$ and any integer $n\geq 0$,
the complex $H\otimes T^{\otimes n}$ is acyclic;
\item[(b)] for any object $F$ in $\G[T]$, there exists a
quasi-isomorphism $H\To F$, with $H$ a $\G$-cofibrant complex,
such that for any object $F'$ in $\G[T]$, the map
$H\otimes F'\To F\otimes F'$ is a quasi-isomorphism.
\end{itemize}

A quasi-isomorphism $H\To F$ satisfying the property
described in (b) will be called a \emph{weakly flat resolution
of $F$}.

Suppose that $T$ is weakly flat.
For each object $F$ in $\G[T]$, choose a weakly flat resolution
$u_{F}: H_{F}\To F$ satisfying the condition (b) above.
Define $\H[T]$ to be $\H\cup\H'$, where $\H'$ is the set
of complexes of shape $\mathrm{Cone}(u_{F})$.
\end{num}

\begin{prop}\label{enlargeweaklyflatdescent}
The couple $(\G[T],\H[T])$ is a weakly flat descent structure.
\end{prop}

\begin{proof}
Note that the elements of $\H[T]$ are $\G[T]$-cofibrant
and acyclic. To prove we have defined a descent structure, it
remains to prove that any $\H[T]$-flasque complex is $\G[T]$-local.
Let $C$ be an $\H[T]$-flasque complex. It is in particular $\H$-flasque,
so that for any $\G$-cofibrant complex $H$ (e.g. $H=H_{F}$),
we have
$$\Hom_{\Htp(\A)}(H,C)=\Hom_{\Der(\A)}(H,C).$$
We also have distinguished triangles
$$H_{F}\To F\To\mathrm{Cone}(u_{F})\To H_{F}[1]$$
in the homotopy category $\Htp(\A)$, so that we deduce
that $C$ is $\G[T]$-local.

To prove the weak flatness, we first notice that
condition (a) implies that for any $n\geq 0$,
the functor $C\longmapsto C\otimes T^{\otimes n}$
is a left Quillen functor from the $\G$-model structure
to the injective model structure.
Condition (b) completes the proof.
\end{proof}

\begin{ex}\label{exweaklyflatobject}
Let $T$ be an object of $\A$. Suppose that
for any integer $n\geq 0$, there exists a short exact
sequence of the form
$$0\To A_{n}\To B_{n}\To T^{\otimes n}\To 0,$$
with $A_{n}$ and $B_{n}$ in $\G$, and
such that for any object $F$ in $\G[T]$, the sequence
$$0\To A_{n}\otimes F\To B_{n}\otimes F\To T^{\otimes n}\otimes F\To 0$$
remains exact. Then $T$ is weakly flat.

Condition (a) is easy to check: if $H$ is a $\G$-cofibrant complex,
then it is degreewise a direct factor of a sum of
elements of $\G$, so that we get short exact sequences
$$0\To A_{n}\otimes H\To B_{n}\otimes H\To T^{\otimes n}\otimes H\To 0.$$
If moreover $H$ is acyclic, as the complexes $A_{n}\otimes H$ and $B_{n}\otimes H$
are then acyclic, the complex $T^{\otimes n}\otimes H$
is acyclic as well. It remains to verify condition (b).
For an object $E$ in $\G$ and an integer $n$, define $H$ to be the complex
$$H=\mathrm{Cone}\big(A_{n}\otimes E\To B_{n}\otimes E\big).$$
The obvious quasi-isomorphism $H\To T^{\otimes n}\otimes E$
is a weakly flat resolution of $F=T^{\otimes n}\otimes E$.
\end{ex}

\section{Localization of derived categories}

\begin{num}\label{cofreplacementT}
From now on, we consider an abelian Grothendieck category $\A$ endowed with
a descent structure $(\G,\H)$.\\
\indent Let $\T$ be an essentially small set of complexes over $\A$.
For any $A$ in $\T$, we choose a quasi-isomorphism $T'\To T$ with $T'$
$\G$-cofibrant (this exists using the cofibrant resolutions of the $\G$-model structure on $\A$), and set
$$\T'=\{T' \ | \ T \in \T\} \ .$$
A morphism of complexes $X\To Y$ over $\A$ is a \emph{$\T$-\'equivalence}
if for any $\H\cup\T'$-flasque complex $K$, the map
$$\Hom_{\Der(\A)}(Y,K)\To\Hom_{\Der(\A)}(X,K)$$
is bijective.
\end{num}

\begin{lm}\label{carTflasque}
Let $K$ be complex over $\A$. The following conditions are equivalent.
\begin{itemize}
\item[(i)] The complex $K$ is $\H\cup\T'$-flasque.
\item[(ii)] The complex $K$ is $\G$-local, and for any $T$ in $\T$, any integer $n$, the group $\Hom_{\Der(\A)}(T,K[n])$
vanishes.
\item[(iii)] The complex $K$ is $\G$-local, and for any $\T$-equivalence $X\To Y$, the map
$$\Hom_{\Der(\A)}(Y,K)\To\Hom_{\Der(\A)}(X,K)$$
is bijective.
\end{itemize}
\end{lm}

\begin{proof}
We know that $K$ is $\H$-flasque if and only if it is $\G$-local, and it is clear
that $K$ is $\H\cup\T'$-flasque if and only if it is $\H$-flasque and $\T'$-flasque.
Moreover, if $K$ is $\G$-local, for any quasi-isomorphism $T'\To T$
with $T'$ $\G$-cofibrant, the map
$$\Hom_{\Htp(\A)}(T,K)\To\Hom_{\Htp(\A)}(T',K)$$
is bijective. This implies that the conditions (i) and (ii) are equivalent.
As (i) implies (iii) by definition, it remains to show that (iii) implies (i).
But this comes from the obvious fact that for any $T'$ in $\T'$, the
map $T'\To 0$ is a $\T$-equivalence.
\end{proof}

\begin{prop}\label{carTequDA}
Let $K$ be a complex over $\A$. The following conditions are equivalent.
\begin{itemize}
\item[(i)] For any $T$ in $\T$, and any integer $n$, the group $\Hom_{\Der(\A)}(T,K[n])$
vanishes.
\item[(ii)] For any $\T$-equivalence $X\To Y$, the map
$$\Hom_{\Der(\A)}(Y,K)\To\Hom_{\Der(\A)}(X,K)$$
is bijective.
\end{itemize}
\end{prop}

\begin{proof}
Let $K\To K'$ be a quasi-isomorphism with $K'$ $\G$-local -- 
this exists using the fibrant resolutions of the $\G$-model structure. 
It is clear that condition (i) (resp. (ii))
is verified for $K$ if and only if it is for $K'$. 
So we can assume that $K$ is $\G$-local.
This proposition then follows immediately from lemma \ref{carTflasque}.
\end{proof}

\begin{paragr}
We will say that a complex $K$ is \emph{$\T$-local} if it satisfies
one of the equivalent conditions of the proposition above.
It is straightforward to check that a map of complexes
$X\To Y$ is a $\T$-equivalence if and only if for any
$\T$-local complex $K$, the map $\Hom_{\Der(\A)}(Y,K)\To\Hom_{\Der(\A)}(X,K)$
is bijective. Hence the notion of $\T$-equivalence
doesn't depend on the descent structure $(\G,\H)$.
\end{paragr}

\begin{prop}\label{propdefcmfloc}
The category $\Comp(\A)$ admits a proper cellular
model category structure with $\T$-equivalences as weak equivalences and
$\G$-cofibrations as cofibratons. This model category structure on $\Comp(\A)$
is the left Bousfield localization of the $\G$-model structure by the maps $0\To T[n]$
for $T$ in $\T$ and $n$ in $\mathbb{Z}$. We will call this model
category structure the \emph{$\G$-model structure associated to $\T$}.\\
\indent The fibrant objects for this model
structure can be characterized as the $\G$-local and $\T$-local
complexes.\\
\indent Let $\Der^{}_\T(\A)$ denotes the localization of $\Comp(\A)$
by the $\T$-equivalences. If $\TDer(\A)$ is the localizing subcategory of $\Der(\A)$
generated by the objects $T[n]$ for $T$ in $\T$ and any integer $n$\footnote{That is
the smallest full triangulated subcategory of $\Der(\A)$ stable by direct sums that contains $\T$.},
then $\Der^{}_\T(\A)=\Der(\A)/\TDer(\A)$. In particular the category $\Der^{}_\T(\A)$
is triangulated. Moreover, the localization functor
$$\Der(\A)\To\Der^{}_\T(\A)$$
is triangulated, and it has a right adjoint that is fully faithful
and identifies $\Der^{}_\T(\A)$ with the full subcategory
of $\Der(\A)$ that consists of $\T$-local complexes.
\end{prop}

\begin{rem}\label{remcarfibgloc2}
It is possible to give an explicit description of the
fibrations for this model category structure; see Corollary \ref{characfibrations3}.
\end{rem}

\begin{proof}
As the $\G$-model structure is (left) proper and cellular, 
we can apply Hirsch\-horn's general construction~\cite{Hir} 
and obtain a cellular left proper model category structure on $\Comp(\A)$ 
as the left Bousfield localization
of the $\G$-model structure 
by the maps $0\To T[n]$ for $T$ in $\T$ and $n$ in $\mathbb{Z}$.
By definition, the cofibrations of this latter model structure still are the $\G$-cofibrations.
There remains to show three points to conclude: characterize the fibrant objects
as annonced in the statement of this proposition, identify $\Der^{}_\T(\A)$
with a full subcategory of $\Der(\A)$, and prove that this model
category structure is right proper. Let us show the first point.
For any pair $(X,Y)$ of objects of $\Comp(\A)$, we can define functorially
(at least in the homotopy category of simplicial sets) a pointed simplicial set
$\derR\Hom(X,Y)$ whose higher homotopy groups are the
$\Hom_{\Der(\A)}(X[n],Y)$'s for $n\geq 0$. In Hirschhorn's construction,
the fibrant objects of the Bousfield localization of the $\G$-model structure are
then characterized as the objects $K$ such that for any integer $n$,
$$\derR\Hom(T[n],K)\To\derR\Hom(0[n],K)\simeq 0$$
is a (weak) homotopy equivalence. But this is is clearly equivalent to say that
the groups $\Hom_{\Der(\A)}(T,K[n])$ vanish. The general theory of left Bousfield localization
gives that the localization
functor $\Der(\A)\To\Der^{}_\T(\A)$ has a fully faithful right adjoint
whose essential image is made of fibrant objects, that are here
exactly the $\G$-local $\T'$-flasque complexes. To identify
$\Der^{}_\T(\A)$ with the quotient $\Der(\A)/\TDer(\A)$, by proposition \ref{carTequDA}
and the theory of left Bousfield localization of triangulated categories (see \cite{Nee}),
it is sufficient to prove that for any $C$ in $\TDer(A)$, the map $0\To C$ is a $\T$-equivalence.
This is an easy exercise left to the reader. To prove right properness, consider now
a pullback square of complexes
$$\xymatrix{
X'\ar[r]^u\ar[d]_{p'}&X\ar[d]^p\\
Y'\ar[r]_v&Y}$$
where $p$ is a fibration for the localized model structure, and $v$ is a $\T$-equivalence.
We have to prove that $u$ is a $\T$-equivalence. But $p$ is also a fibration
for the $\G$-model structure, and as the $\G$-model structure is right proper,
this square is homotopy cartesian with respect to the $\G$-model structure.
This is equivalent to say that we have a canonical distinguished square
$$X'\To Y'\oplus X\To Y\To X'[1]$$
in $\Der(\A)$. As the localization functor from $\Der(\A)$ to $\Der^{}_\T(\A)$ is triangulated,
this is also a distinguished square in $\Der^{}_\T(\A)$. But $v$ is an isomorphism in
$\Der^{}_\T(\A)$, which then implies that $u$ is also an isomorphism.
As a morphism in a model category is a weak equivalence if and only if it is an
isomorphism in the corresponding homotopy category, this implies
that $u$ is a $\T$-equivalence.
\end{proof}

\begin{prop}\label{carrehtpcartTlc}
A pushout square of complexes of $\A$
$$\xymatrix{
K\ar[r]^{u}\ar[d]_{i}&L\ar[d]^j\\
K'\ar[r]_{u'}&L'
}$$
in which $i$ is a monomorphism is a homotopy pushout square
with respect to the $\G$-model structure associated to $\T$.
In particular, if $i$ (resp. $u$) is a $\T$-equivalence, then
$j$ (resp. $u'$) is a $\T$-equivalence.
\end{prop}

\begin{proof}
As the identity of $\Comp(\A)$ is a left Quillen
functor from the $\G$-model structure to the
$\G$-model structure associated to $\T$,
it is sufficient to prove that such a pushout square
is a homotopy pushout square with respect to the
$\G$-model structure. But the notion of homotopy
pushout square only depends on the weak equivalences,
so that we are free to work with the injective model structure.
The assertion is then trivial.
\end{proof}

\begin{prop}\label{carrehtpcartTlc2}
The $\T$-equivalences are stable by filtering
colimits in $\Comp(\A)$.
\end{prop}

\begin{proof}
We know that the quasi-isomorphisms
are stable by filtering colimits (just because the
filtering colimits are exact in any Grothendieck
abelien category). This is equivalent to say that
for any filtering category $I$ and any functor $F$
from $I$ to $\Comp(\A)$, the canonical map
$$\hocolim_{i}F_{i}\To\varinjlim_{i}F_{i}$$
is a quasi-isomorphism. But as the functor
$\hocolim$ sends termwise $\T$-equivalences to
$\T$-equivalences, this implies our assertion.
\end{proof}

\begin{prop}\label{fibrationslocalisees}
A morphism of complexes of $\A$ is a fibration with respect
to the $\G$-model structure associated to $\T$
if and only if it is a fibration with respect to the
$\G$-model structure with $\T$-local kernel.
\end{prop}

\begin{proof}
Let $p:K\To L$ be a fibration with respect
to the $\G$-model structure associated to $\T$.
It has to be a fibration with respect to the
$\G$-model structure and as the fibrations are stable by pullback,
the map from the kernel of $p$ to $0$ has to be
a fibration with respect
to the $\G$-model structure associated to $\T$.
Hence the kernel of $p$ is $\T$-local.
Conversely, let $p$ be a fibration with respect
to the $\G$-model structure with $\T$-local kernel.
We can factor $p=qi$ where $i:K\To M$ is
a $\G$-cofibration and a $\T$-equivalence and
$q:M\To L$ is a fibration with respect
to the $\G$-model structure associated to $\T$.
As the $\G$-model structure associated to $\T$
is right proper, the induced map from the kernel
of $p$ to the kernel of $q$ is a $\T$-equivalence.
But as these two kernels are $\T$-local, the map
$\mathrm{ker}\, p \To \mathrm{ker}\, q$ has to be a
quasi-isomorphism. But both $p$ and $q$
are epimorphisms (this comes
from the fact that they have the right lifting property with respect to
the maps $0\To D^nE$ for $E$ in $\G$ and $n$ an integer).
We thus have a morphism of distinguished triangles in
$\Der(\A)$
$$\xymatrix{
{\mathrm{ker}}\, p\ar[r]\ar[d]_{j}&K\ar[r]^p\ar[d]_{i}&L\ar[r]\ar[d]^{1_{L}}
&{\mathrm{ker}}\, p[1]\ar[d]^{j[1]}\\
{\mathrm{ker}}\, q\ar[r]&M\ar[r]_{q}&L\ar[r]&{\mathrm{ker}}\, q[1]}$$
where $j$ is an quasi-isomorphism. This proves that $i$ is
also a quasi-isomorphism. As it is also a $\G$-cofibration,
it follows that $p$ is a retract of $q$. As $q$ is a fibration with respect
to the $\G$-model structure associated to $\T$, $p$ must
have the same property.
\end{proof}

We recall here a folklore result that is quite useful to play the game
between Bousfield localizations and Quillen functors; see \cite{Hir}.

\begin{lm}\label{1locbousfieldcmf}
Let $G:\C\To \C'$ be a left Quillen functor between
model categories. Suppose that we have a class $S$
of maps in $\C$ such that the left Bousfield
localization of $\C$ with respect to $S$ exists. 
We denote by $\C_{S}$
the category $\C$ endowed with the model category structure
obtained as the left Bousfield localization by $S$.
Denote by $D$ a right adjoint to $G$.
Then the following conditions are equivalent.
\begin{itemize}
\item[(a)] The functor $G$ is a left Quillen functor from $\C_{S}$ to $\C'$.
\item[(b)] For any fibrant object $X'$ in $\C'$, $D(X')$
is a fibrant object in $\C_{S}$.
\item[(c)] The functor $G$ sends the weak equivalences between
cofibrant objects in $\C_{S}$ to weak equivalences.
\item[(d)] The functor $\derL G : \ho(\C)\To\ho(\C')$
sends the elements of $S$ to isomorphisms.
\end{itemize}
\end{lm}

\begin{prop}\label{fonct2}
Let $\A$ (resp. $\A'$) be a Grothendieck abelian category
 with descent structure $(\G,\H)$ (resp. $(\G',\H)$)
Let $f^*:\A'\To \A$ be a functor which has a right adjoint $f_*$
and satisfies descent.
Consider a set $\T$ of $\G$-cofibrant complexes over $\A$, and a set $\T'$
of $\G'$-cofibrant complexes over $\A'$
such that $f^*(\T')\subset\T$.
Then the pair of adjoint functors
$$f^*:\Comp(\A')\To\Comp(\A)\quad\text{and}
\quad f_*:\Comp(\A)\To\Comp(\A')$$
is a Quillen pair with respect to the $\G'$-model structure
associated to $\T'$ and the $\G$-model structure associated to $\T$.
In particular, we have a pair of adjoint functors
$$\derL f^*:\Der_{\T'}(\A')\To\Der_{\T}(\A)\quad\text{and}\quad
\derR f_*:\Der_{\T}(\A)\To\Der_{\T'}(\A')$$
obtained as the left and right derived functors
of $f^*$ and $f_*$ respectively with respect
to the above model structures.
\end{prop}

\begin{proof}
By virtue of theorem \ref{fonct1}, the functor $f^{*}$ is a left Quillen functor,
so that it has a total left derived functor
$$\derL f^{*}:\Der(\A')\To\Der_{\T}(\A) \ .$$
Using Proposition \ref{propdefcmfloc} and
the condition (d) of Lemma \ref{1locbousfieldcmf},
it is sufficient to check that $\derL f^{*}$ sends the elements
of $\T'$ to null objects in $\Der_{\T}(\A)$, which is obvious.
\end{proof}

\begin{rem}
We have of course the same 2-functoriality
results as in remark \ref{composefctder}
in the case the obvious compatibilities
are verified.
\end{rem}

\begin{num}\label{defflatgcofcomploc}
Let $\A$ be a closed symmetric monoidal
Grothendieck category endowed with a
weakly flat descent structure $(\G,\H)$.
A set $\T$ of $\G$-cofibrant complexes of $\A$
will be said to be \emph{flat} (with respect to $\G$)
if for any $E$ in $\G$ and any $T$ in $\T$, $E\otimes T$ is in $\T$.
\end{num}

\begin{cor}\label{cmf4}
Let $\A$ be a closed symmetric monoidal Grothendieck category.
We suppose that we have a weakly flat descent structure $(\G,\H)$
on $\A$ and a flat set $\T$ of $\G$-cofibrant complexes of $\A$.
Then the $\G$-model structure associated to $\T$ is
a symmetric monoidal model category, and
the localization functor $\Der(\A)\To\Der_{\T}(\A)$ is
a symmetric monoidal functor.

If moreover $\G$ is flat, then this model category satisfies the
monoid axiom of \cite{SS}. 
\end{cor}

\begin{proof}
The proof of Lemma \ref{cmf33} can be modified slightly to
prove that for any $\G$-cofibrant complex $A$ and any $T$ in $\T$,
the map $0 \To A\otimes T$ is a $\T$-equivalence.
The first assertion thus follows from Proposition \ref{carflat1}
and from a new application of Lemma \ref{1locbousfieldcmf}.
If we assume that $\G$ is flat, then $(\G,\H)$ is flat (\ref{cmf333}).
To prove the second assertion, we will consider the left Bousfield
localization of the injective model structure
on the category $\Comp{(\A)}$ by the maps
$0\To T[n]$ for any integer $n$ and any $T$ in $\T$.
We will call the latter model category structure
the $\T$-injective model structure.
The identity of $\Comp{(\A)}$ is a left Quillen
equivalence from the $\G$-model structure associated
to $\T$ to the $\T$-injective model structure.
One deduces easily from this, from Propositions \ref{carrehtpcartTlc}
and \ref{carrehtpcartTlc2}, and from the identification
$\Der^{}_\T(\A)=\Der(\A)/\TDer(\A)$ of Proposition \ref{propdefcmfloc}
that the class of weak equivalences coincides for the two
model structures associated to $\T$. In particular,
the class of monomorphisms which are
$\T$-equivalences is stable by transfinite
composition and by pushout. It is then
sufficient to prove that for any complex $C$,
the functor $C\longmapsto C\otimes D$
is a left Quillen functor from the $\G$-model structure associated to $\T$
to the $\T$-injective model structure. It follows from
Proposition \ref{monoidaxiom} and from lemma \ref{1locbousfieldcmf}
that it is sufficient for this to prove that for any $\T$-equivalence
$A\To B$ between $\G$-cofibrant complexes, the map
$C\otimes A\To C\otimes B$ is a $\T$-equivalence.
But for any $\G$-cofibrant complex $K$, $C\otimes K$
is isomorphic to the derived tensor product of $C$ and $K$
(see Remark \ref{tensdercoincide}), so that by functoriality
$C\otimes A\To C\otimes B$ is an isomorphism in
the homotopy category $\Der^{}_\T(\A)$. Hence the result.
\end{proof}

\begin{ex} \label{DMeff}
Consider the notations of Example \ref{exshtrans} and \ref{Dftr}.

For a regular scheme $S$, we let $\T_S$ be the essentially small set 
of complexes of sheaves with transfers of the form
$L_S(\AA^1_X) \xrightarrow{p_*} L_S(X)$
for any smooth $S$-scheme $X$,
 $p:\AA^1_X \rightarrow X$ being the canonical projection.

We define the category of motivic complexes over $S$
as the category $$\DMe(S):=\Der_{\T_S}(\ftr S)$$ with the
notation of Proposition \ref{propdefcmfloc}.

A complex $C$ is $\T_S$-local, simply called $\AA^1$-local,
if its Nisnevich hypercohomology sheaves are homotopy invariant.
Remark that if $C$ is Nisnevich fibrant,
then it is $\AA^1$-local if for any smooth $S$-scheme $X$,
the map induced by the projection $C(X) \rightarrow C(\AA^1_X)$ 
is a quasi-isomorphism, or equivalently, the cohomology presheaves
of $C$ are homotopy invariant.

From Proposition \ref{propdefcmfloc}, the category $\DMe(S)$
is equivalent to the full subcategory of $\Der(\ftr S)$
made by the Nisnevich fibrant and $\AA^1$-local complexes.
This result has to be compared with the original definition
of Voevodsky of motivic complexes over a perfect field
(see \cite{FSV}, chapter 5, following proposition 3.1.13).
In fact, according to the strongest result of the theory
(theorem 3.1.12 of \emph{loc. cit.}), over a perfect field,
a (bounded below) complex is $\AA^1$-local if its cohomology
sheaves are homotopy invariant\footnote{Among the notable facts
this theorem implies is the description of the homotopy $t$-structure
on $\DMe(k)$ and the Gersten resolution for an object
of its heart}.

According to Corollary \ref{cmf4} and Example \ref{DMmonoidal},
the category $\DMe(S)$ is symmetric monoidal.
According to Proposition \ref{fonct2} and Example \ref{DM_4functors},
for any morphism $f:T \rightarrow S$ of regular schemes,
we have a pair of adjoint functors
$$
\derL f^*:\DMe(S) \leftrightarrows \DMe(T):\derR f_*
$$
and if $f$ is smooth, another pair
$$
\derL f_\sharp:\DMe(T) \leftrightarrows \DMe(S):\derR f^*=f^*.
$$
If we denote by $M_S(X)$ the object represented by a smooth scheme
in $\DMe(S)$, we have the fundamental relations
$$
M_S(X) \otimes^\derL M_S(Y)=M_S(X \times_S Y), \
 \derR f^* M_S(X)=M_T(X \times_S T), \
  \derL f_\sharp M_T(Y)=M_S(Y)
$$
whenever it makes sense.
\end{ex}

\section{Presentations of the derived category}

\begin{paragr}
Let $R$ be a commutative ring with unit, and $\AA$
an essentially small additive $R$-linear category.
We denote by $\AA\Mod$
the category of right $\AA$-modules, that is the category
of $R$-linear additive functors from $\op{\AA}$
to the category of $R$-modules. The Yoneda embedding of
$\AA$ in $\AA\Mod$ is additive, $R$-linear and fully faithful.
For any object $X$ of $\AA$, we will also denote by $X$
the right $\AA$-module represented by $X$.
\end{paragr}

\begin{prop}\label{basiccmfmodules}
Let $\G_{\AA}$ be a set of representative objects of $\AA$.
Then $(\G_{\AA},\varnothing)$ is a descent structure
on the category of $\AA$-modules. In other words,
the $\G_{\AA}$-model structure on the category of
complexes of $\AA$-modules is the model category structure
with the quasi-isomorphisms as weak equivalences and the
epimorphisms as fibrations.
\end{prop}

\begin{proof}
It is clear that $(\G_{\AA},\varnothing)$ is a descent structure
on the category of $\AA$-modules because $\G_{\AA}$
is obviously a set of generators and because any
complex of $\AA$-modules is $\G_{\AA}$-local.
By definition of the $\G_{\AA}$-model structure,
the fibrations are the maps which have
the right lifting property with respect to the maps
$0\To D^nE$ for any $E$ in $\G_{\AA}$ and any integer $n$.
But as $D^nE$ represents the evaluation
functor $K\longmapsto\Hom(E,K_{n})$, this implies
that the fibrations are exactly the epimorphisms
of complexes of $\AA$-modules.
\end{proof}

\begin{paragr}\label{defAAKanext}
Let $\A$ be a Grothendieck category and $(\G,\H)$
a descent structure on $\A$. Let $\AA$
be the full subcategory of $\A$ whose objects are the
finite direct sums of elements of $\G$, and let $i$ be the
inclusion functor. We have a functor
$$i^* : \A\To\AA\Mod$$
defined by
$$i^*(F)(X)=\Hom_\A(i(X),F)$$
for any object $F$ of $\A$ and any object $X$ of $\AA$.
The functor $i^*$ has a left adjoint $i_!$
that can be defined as the left Kan extension of $i$.
$$i_! : \AA\Mod\To\A$$
In particular, if $X$ is an object of $\AA$,
we have $i_!(X)=i(X)$ (recall that we denote also by $X$
the presheaf on $\AA$ represented by $X$).\\
\end{paragr}

\begin{lm}\label{quillenpref}
The functor $i_{!}:\Comp(\AA\Mod)\To\Comp(\A)$
is a left Quillen functor from the $\G_{\AA}$-model
structure (\ref{basiccmfmodules}) to the $\G$-model structure.
\end{lm}

\begin{proof}
As $i_{!}$ preserves colimits, it is sufficient to prove that
it sends the generating cofibrations (resp. trivial cofibrations)
to cofibrations (resp. to trivial cofibrations). But for any object
$X$ of $\AA$, $E=i_{!}(X)=i(X)$ is in $\G$ and for any integer $n$
we have $i_{!}(D^nX)=D^nE$ and $i_{!}(S^nX)=S^nE$, so that
conclusion follows.
\end{proof}

\begin{lm}\label{equivGcofibrants}
The functor $i_{!}:\Comp(\AA\Mod)\To\Comp(\A)$
induces an equivalence of categories between the
category of $\G_{\AA}$-cofibrant complexes of $\AA$-modules
to the category of $\G$-cofibrant complexes of $\A$.
\end{lm}

\begin{proof}
The functor $i_{!}$ is fully faithful on the full subategory
of $\Comp(\AA\Mod)$ of the $\G_{\AA}$-cofibrant complexes.
This comes from the facts that $i$ is fully faithful by definition
and that if $C$ is a $\G_{\AA}$-cofibrant
complex of $\AA$-modules, then for any integer $n$, $C^n$
is a direct summand of direct sum of representable presheaves
on $\AA$.
The problem is thus to describe the essential image of $i_{!}$
restricted to $\G_{\AA}$-cofibrant complexes.
This will be proved inductively using the
small object argument. For this purpose, we have to introduce
the following class of complexes. Let $\C$ be the class of complexes $K$ of $\A$
such that there exists  a $\G_{\AA}$-cofibrant complex
of $\AA$-modules $L$ and an isomorphism $i_{!}(L)\simeq K$.
We will show that the class $C$ is stable by all the operations
we need for the small object argument.
The class $\C$ contains all the elements of $\G$,
and is stable by direct sums. We deduce from the
fact that  $i_{!}$ is fully faithful on
$\G_{\AA}$-cofibrant complexes  that
any direct factor of an object in $\C$ is in $\C$.\\
\indent Say that a morphism $u:K\To L$ of complexes of $\A$
is an \emph{elementary $\G$-cofibration}
if there exists a pushout square of shape
$$\xymatrix{
{\bigoplus_{j\in J}S^{n_{j}+1}E_{j}}\ar[r]\ar[d]&K\ar[d]\\
{\bigoplus_{j\in J}D^{n_{j}}E_{j}}\ar[r]&L}$$
where for any $j$ in $J$, $E_{j}$ is in $\G$ and $n_{j}$
is an integer. For such a map $u:K\To L$, if moreover $K$
is isomorphic to some $i_{!}(A)$ with a $\G_{\AA}$-cofibrant complex $A$,
and if one sets $X_{i}=i^*(E_{i})$,
we get a pushout square in the category of complexes
of $\AA$-modules
$$\xymatrix{
{\bigoplus_{j\in J}S^{n_{j}+1}X_{j}}\ar[r]\ar[d]&i^*(K)\ar[d]^v\\
{\bigoplus_{j\in J}D^{n_{j}}X_{j}}\ar[r]&M.}$$
One checks easily that $i_{!}(M)\simeq L$ and that
$i^*(u)\simeq v$. Hence $L$ is in $\C$ and
$i^*(u)$ is a $\G_{\AA}$-cofibration. In conclusion, we have proved that
if $u:K\To L$ is an elementary cofibration with $K$ in $\C$,
there exists a (unique) $\G_{\AA}$-cofibration between
$\G_{\AA}$-cofibrant objects $v$ such that $u\simeq i_{!}(v)$.
Let now $\lambda$ be an ordinal with initial element $0$,
and consider a colimit preserving functor $F : \lambda\To\Comp(\A)$
(that is a $\lambda$ sequence as defined in \cite[2.1.1]{Hov})
such that for any $\gamma<\lambda$, $F(\gamma)$ is in $\C$
and $F(\gamma)\To F(\gamma+1)$ is an elementary cofibration.
Then one can find a colimit preserving functor $G:\lambda\To\Comp(\AA\Mod)$
and an isomorphism $i_{!}(G)\simeq F$ (this is obvious when $\lambda$
is finite, and we do it by transfinite induction in the general case).
But then the colimit of $G$ is $\G_{\AA}$-cofibrant, and we get
$$\varinjlim_{\gamma<\lambda}F({\gamma})\simeq
\varinjlim_{\gamma<\lambda}i_{!}\big(G({\gamma})\big)\simeq
i_{!}\Big(\varinjlim_{\gamma<\lambda}(G({\gamma})\Big) \ . $$
In other words, $\smash{\varinjlim_{\gamma<\lambda}F({\gamma})}$
is in $\C$. As $\C$ contains all the elements of $\G$, the small
object argument applied to the set of generating $\G$-cofibrations implies that
$\C$ contains all the $\G$-cofibrant complexes. 
\end{proof}

\begin{paragr}
We consider the set $\H_{\AA}=\{i^*(H) \ , \, H \in \H \}$.
It follows from the preceding Lemma that any element of $\H_{\AA}$
is $\G_{\AA}$-cofibrant. Moreover, for any $H$, the canonical map $i_{!}i^*(H)\To H$
is an isomorphism of complexes.
\end{paragr}

\begin{thm}\label{quilenequivpref}
The functor $i_{!}:\Comp(\AA\Mod)\To\Comp(\A)$ is a Quillen equivalence from
the $\G_{\AA}$-model structure associated to $\H_{\AA}$ on $\Comp(\AA\Mod)$
to the $\G$-model structure on $\Comp(\A)$.
\end{thm}

\begin{proof}
Proposition \ref{fonct2} shows that we have defined a left Quillen functor.
It remains to prove that the total left derived functor
$$\derL i_{!}: \Der_{\H_{\AA}}(\AA\Mod)\To\Der(\A)$$
is an equivalence of categories.
As any complex of $\A$ is isomorphic in $\Der(\A)$ to a $\G$-cofibrant
complex, the essential surjectivity is induced by Lemma \ref{equivGcofibrants}.
Let $K$ and $L$ be two complexes of $\AA$-modules.
We want to show that the map
$$\Hom_{\Der_{\H_{\AA}}(\AA\Mod)}(K,L)\To\Hom_{\Der(\A)}(\derL i_{!}(K),\derL i_{!}(L))$$
is an isomorphism of abelian groups. For this purpose, we can suppose that
both $K$ and $L$ are $\G_{\AA}$-cofibrant and that $L$
is fibrant with respect to the $\G_{\AA}$-model structure associated to $\H_{\AA}$.
We can prove that $i_{!}(L)$ is $\G$-local: this follows from the fact
that $i_{!}$ is fully faithful on $\G_{\AA}$-cofibrant complexes
and from the fact that $L$ is $\H_{\AA}$-local
(or equivalently $\H_{\AA}$-flasque), which implies
that $i_{!}(L)$ is $\H$-flasque. Whence $i_{!}(L)$ is $\G$-local
by Theorem \ref{cmf2}. But we have the following canonical identifications.
$$\begin{aligned}
\Hom_{\Der_{\H_{\AA}}(\AA\Mod)}(K,L)&=\Hom_{\Htp(\AA\Mod)}(K,L)\\
&=\Hom_{\Htp(\A)}(i_{!}(K),i_{!}(L))\quad\text{(by Lemma \ref{equivGcofibrants})}\\
&=\Hom_{\Der(\A)}(i_{!}(K),i_{!}(L))\\
&=\Hom_{\Der(\A)}(\derL i_{!}(K),\derL i_{!}(L)) \ .
\end{aligned}$$
This achieves the proof.
\end{proof}

\begin{paragr}\label{defGsurjective}
Say that a morphism $u:K\To L$
is a \emph{$\G$-surjection} if for any $E$ in $\G$, the map
$$\Hom_{\A}(E,K)\To\Hom_{\A}(E,L)$$
is surjective.
\end{paragr}

\begin{cor}\label{characfibrations2}
A morphism of complexes of $\A$ is a fibration for the $\G$-model structure
if and only if it is a degreewise $\G$-surjection with
$\G$-local kernel.
\end{cor}

\begin{proof}
Any fibration for the $\G$-model structure on $\Comp(\A)$
is a $\G$-surjection with $\G$-local kernel.
Conversely, it follows from Theorem \ref{cmf2} and Propositions
\ref{fibrationslocalisees} and \ref{basiccmfmodules}
that a morphism $p:X\To Y$ of complexes of $\A$ is degreewise
$\G$-surjective with $\G$-local kernel if and only if
$i^*(p)$ is a fibration with respect to the $\G_{\AA}$-model
structure assocated to $\H_{\AA}$. This implies
by Lemma \ref{equivGcofibrants} and Theorem \ref{quilenequivpref}
that $p$ has the right lifting property with respect to
trivial cofibrations between cofibrant objects.
But it follows from Lemma \ref{Icof} that
all the generating trivial cofibrations are cofibrations between
cofibrant objects. Hence the result.
\end{proof}

\begin{cor}\label{characfibrations3}
Let $\A$ be a Grothendieck category endowed with a descent
structure $(\G,\H)$ and a set $\T$ of complexes of $\A$.
A morphism of complexes of $\A$ is a fibration with respect to the
$\G$-model structure associated to $\T$ if and only if it is
degreewise $\G$-surjective with a $\G$-local and $\T$-local kernel.
\end{cor}

\begin{proof}
This follows from Proposition \ref{fibrationslocalisees}
and Corollary \ref{characfibrations2}.
\end{proof}

\begin{rem}\label{generators}
The preceding result implies that
we can describe a generating set of trivial
cofibration of the $\G$-model structure associated to $\T$ on $\Comp(\A)$: just add to the
generating trivial cofibrations of the $\G$-model structure on $\Comp(\A)$
the maps $T'\oplus T'[n]\To \mathrm{Cyl}(T')[n]$
where $n$ runs over the integers, $T$ over the set $\T$, and where
$T'$ is a choice of a $\G$-cofibrant resolution of $T$.
\end{rem}

\section{How to get compact generators}

\begin{num}
Recall that for a triangulated category $\T$, an object $X$ is
\emph{compact} if for any set $I$ and any family of objects
$Y_{i}$, $i\in I$, in $\T$, the map
$$\bigoplus_{i\in I}\Hom_{\T}(X,Y_{i})\To\Hom_{\T}(X,\bigoplus_{i\in I}Y_{i})$$
is surjective (hence bijective). One checks easily that
the full subcategory of compact objects in $\T$ is
a triangulated subcategory of $\T$.

Consider a Grothendieck abelian category $\A$.
We say that a descent structure $(\G,\H)$ on $\A$
is \emph{bounded} if the complexes $H$ in $\H$
are all bounded and degreewise finite sums of
objects in $\G$. The following lemma is
straightforward.
\end{num}

\begin{lm}\label{compactmodulestriv}
Let $R$ be a commutative ring with unit and
$\AA$ be an $R$-linear additive category.
Any bounded complex of $\AA$-modules which is
degreewise a finite sum of representable $\AA$-modules is compact
in the triangulated category  $\Htp(\AA\Mod)$.
\end{lm}

\begin{num} \label{def_Dec_c}
Consider given a bounded descent structure $(\G,\H)$ on $\A$.
Let $\AA$ be the full subcategory of $\A$
that consists of finite sums of objects in $\G$.
We define the triangulated category $\Der_{c}(\A)$
as follows. Consider the bounded homotopy category of $\AA$,
denoted by $\Htp^b(\AA)$. The previous lemma
implies easily that the triangulated
category of compact objects in $\Htp(\AA\Mod)$
contains the pseudo-abelianization\footnote{This
means the idempotent completion: we add formaly
the kernel and cokernel of any projector.} of $\Htp^b(\AA)$.
By virtue of Lemma \ref{equivGcofibrants}, we can consider
the objects in $\H$ as bounded complexes in $\AA$.
Let $<\H>$ be the thick subcategory
of $\Htp^b(\AA)$ generated by the objects in $\H$, and
define $\Der_{c}(\A)$ to be the pseudo-abelianization of the
Verdier quotient of $\Htp^b(\AA)$ by $<\H>$. It follows from a
general result of Balmer and Schlichting~\cite{balsch}
that $\Der_{c}(\A)$ is canonically endowed with
a structure of triangulated category such that the functor
$\Htp^b(\AA)\To\Der_{c}(\A)$ is triangulated.
\end{num}

\begin{thm}\label{presderivedcat}
Any object in $\G$, seen as a complex concentrated in degree $0$,
is a compact object in $\Der(\A)$. Moreover, the functor
$$\derL i_{!}:\Der(\AA\Mod)\To\Der(\A)$$
induces an equivalence of triangulated categories
between the category of compact objects of $\Der(\A)$
and the category $\Der_{c}(\A)$.
\end{thm}

\begin{proof}
Let $\T$ be the localizing subcategory of
$\Htp(\AA\Mod)$ generated by the re\-pre\-sen\-ta\-ble
$\AA$-modules.
It follows easily from Theorem \ref{quilenequivpref}
and from Lemma \ref{equivGcofibrants}
that $\Der(\A)$ is the quotient of $\T$
by the localizing subcategory of $\T$
generated by the objects in $\H$.
But Lemma \ref{compactmodulestriv}
implies that the objects in $\G$ or in $\H$
are compact in $\T$. The proof is thus a direct
consequence of Thomason Theorem; see \cite[Theorem 4.4.9]{Nee}.
\end{proof}

\begin{ex}\label{DMgm}
Consider the notations of Example \ref{DMeff}.

For any regular scheme $S$,
 the descent structure $(\G_S,\H_S)$ is bounded.
Recall the category of finite correspondences $\smc S$
over $S$ from \cite[4.1.19]{Deg7} is the category whose objets
are smooth $S$-schemes and morphisms from $X$ to $Y$
are cycles in $X \times_S Y$ which support is finite
equidimensional over $X$. It is additive and symmetric monoidal.

According to \ref{def_Dec_c}, we consider the thick subcategory
$<\H_S,\T_S>$ of $K^b(\smc S)$ generated by the classes $\H_S$ and 
$\T_S$. We denote by $\DMgme(S)$ the pseudo-abelianization of
$K^b(\smc S)/<\H_S,\T_S>$.

The preceding theorem implies we have a canonical fully faithful
functor $$\DMgme(S) \rightarrow \DMe(S)$$ whose essential image
is the category of compact objects of $\DMe(S)$.
\end{ex}

\section{Stabilization and symmetric spectra}\label{sectionstabilisation}

\begin{paragr}
Let $\A$ be a Grothendieck category. We suppose that $\A$
is endowed with a symmetric monoidal structure. We denote by
$\otimes$ and $\unit$ the tensor product and the unit respectively.
We also suppose that for any object $X$ of $\A$, the functor
$$Y\longmapsto X\otimes Y$$
preserves colimits (in practice, $\A$ will be the category of
(complexes of) sheaves of $R$-modules on a Grothendieck site,
and $\otimes=\otimes_R$ the usual tensor product on such objects).
\end{paragr}

\begin{paragr}
We will consider in this subsection $G$-objects of $\A$ for
various (finite) groups $G$, that are object $A$ of $\A$
equipped with a representation $G \To \mathrm{Aut}_\A(A)$.

We introduce the following notations~:

If $A$ is an object of $\A$, we put
$G \times A=\bigoplus_{g \in G} A$ considered as a $G$-object via the
permutation isomorphisms of the summands.

If $H$ is a subgroup of $G$, and $A$ is an $H$-object, $G \times A$
has two actions of $H$~: the first one, say $\gamma$, is obtained via
the inclusion $H \subset G$, and the second one denoted by $\gamma'$,
is obtained using the structural action of $H$ on $A$, making $H$
acting diagonally on $G \times A$. We define $G \times_H A$ as the
equalizer of the family of morphisms
$(\gamma_\sigma-\gamma'_\sigma)_{\sigma \in H}$, and consider it 
equipped with its induced action of $G$.

We write $\mathfrak S_n$ for the symmetric group (that is the group
of automorphisms of the set with $n$ elements, $n\geq 0$), and
we make the following definition.
\end{paragr}

\begin{df}
A symmetric sequence of $\A$ is a sequence
$(A_n)_{n \in \NN}$ 
such that for each $n \in \NN$, $A_n$ is a $\mathfrak S_n$-object of
$\A$. A morphisms of symmetric sequences is a degreewise
equivariant morphism in $\A$.
\end{df}

\begin{paragr}
For an integer $n\geq 0$, we denote by $\underline n$
the set of integers $i$ such that $1\leq i\leq n$.
We then introduce the category $\mathfrak S$ as follows.
The objects are the sets $\underline n$ for $n\geq 0$, and the
morphisms are the bijections between such sets. Then, $\mathfrak S$ is a
groupoid and the category of symmetric sequences
can be described as the category of functors from $\mathfrak S$
to $\A$. We use the notation $\sseq A$ to
denote the category of symmetric sequences of $\A$. 

We have a full embedding
$$\A \To \sseq A\quad , \qquad X\longmapsto X\{0\}$$
which is a left adjoint to the $0$-evaluation functor
$\sseq A\To \A \, , A_* \longmapsto A_0$. More precisely, if $X$
is an object of $\A$, $X\{0\}$ is the symmetric sequence defined by
$X\{0\}_0=X$ and $X\{0\}_n=0$ for $n>0$.

Consider an integer $i \in \NN$ and a symmetric sequence $A_*$ of $\A$. We put
\begin{equation}\begin{split}
A_*\{-i\}=n \longmapsto \left\{
\begin{array}{ll}
\mathfrak S_n \times_{\mathfrak S_{n-i}} A_{n-i}
 & \text{if } n \geq i \\
0 & \text{otherwise.} 
\end{array}\right.
\end{split}\end{equation}
This defines an endofunctor on $\A^\mathfrak S$, and we have 
$A_*\{-i\}\{-j\} \simeq A_*\{-i-j\}$.
For an object $X$ of $\A$, we will write $X\{-i\}=X\{0\}\{-i\}$.
For $i\geq 0$, the functor $A_*\longmapsto A_*\{-i\}$ has a right adjoint
$A_*\longmapsto A_*\{i\}$
which is defined by $A_*\{i\}_n=A_{n+i}$, where the action of $\mathfrak S_n$
on $A_{n+i}$ is induced by the canonical inclusion of $\mathfrak S_n$
in $\mathfrak S_{n+i}$.

Remark that for any integer $n \in \NN$, the functor
$\A\To \sseq A , X \mapsto X\{-n\}$ is left
adjoint to the $n$-evaluation functor
$\sseq A \To \A, A_* \longmapsto A_n$. Moreover, the
collection of functors
$\sseq A \xrightarrow{Ev_n} \A$
preserves every limits and colimits and is conservative. 

We now define a symmetric monoidal structure on $\sseq A$.

First of all, let us introduce a variant of $\mathfrak S$ by
considering the category $\mathfrak S'$ of finite sets with morphisms
only the bijections. The canonical functor
$\mathfrak S\To\mathfrak S'$ is thus an equivalence of categories. 
This means that the restriction functor
$\A^{\mathfrak S'} \To \A^\mathfrak S$ is an equivalence of categories. 

Let now $E,F:\mathfrak S' \To \A$ be functors. We
define a functor 
\begin{equation}\label{defsymtensoreq}
\begin{split}
\begin{array}{rcl}
E \otimes^\mathfrak S F:\mathfrak S' & \To & \A \\
N & \longmapsto & \bigoplus_{N=P \sqcup Q} E(P) \otimes F(Q)
\end{array}
\end{split}
\end{equation}
This makes $\A^{\mathfrak S'}$, and thus $\sseq A$ into a
symmetric monoidal category, the unit being the symmetric
sequence $\unit\{0\}$, and the symmetry isomorphism being induced by
the symmetry isomorphism of $\otimes$ in $\A$. It is then straightforward
to check that the canonical functor $X\longmapsto X\{0\}$ is a symmetric
monoidal functor from $\A$ to $\sseq A$.
\end{paragr}

\rem The above definition is only a way to simplify the combinatorics
of $\otimes^\mathfrak S$. To help the reader using it, we point out
that to give a morphism of symmetric objects
$A_* \otimes^\mathfrak S B_* \To C_*$ is equivalent to give
families of maps $A_p \otimes B_q \To C_{p+q}$ which are
$\mathfrak S_p \times \mathfrak S_q$-equivariant, the action on the
right hand side being given by the canonical inclusion
$\mathfrak S_p \times \mathfrak S_q \To \mathfrak S_{p+q}$. 

\begin{paragr}\label{defunivsymmmnoid}
Let $S$ be an object of $\A$. We define $S^{\otimes n}$ as
$$S^{\otimes n}=
\underset{{\text{$n$ times}}}{{\underbrace{S\otimes \dots \otimes S}}} \ .$$
The group of permutations $\mathfrak S_n$ acts on $S^{\otimes n}$ by the
structural symmetry iso\-mor\-phisms of the
monoidal structure. We denote by $Sym(S)$ the symmetric sequence
$(S^{\otimes n})_{n \in \NN}$ obtained by considering the above
actions.

Moreover, the canonical isomorphism functor
$$
S^{\otimes p} \otimes S^{\otimes q} \To S^{\otimes p+q}
$$
obtained using the associativity isomorphism, is $\mathfrak S_p \times
\mathfrak S_q$-equivariant. Thus, this defines a morphism
$$
Sym(S) \otimes^\mathfrak S Sym(S) \To Sym(S).
$$
One checks now that with this structural morphism, $Sym(S)$ is a
commutative monoid in the category $\A^\mathfrak S$ (the unit
is defined by the identity of $\unit$).
Recall the symmetric sequence $S\{0\}$ defined by
$S\{0\}_0=S$ and $S\{0\}_n=0$ for $n>0$, and remember that
$S\{-1\}=S\{0\}\{-1\}$.
So $S\{-1\}$ is the symmetric sequence defined by
$S\{-1\}_n=0$ if $n\neq 1$ and $S\{-1\}_1=S$.
One has a natural morphism of symmetric sequences
$$S\{-1\}\To Sym(S)$$
induced by the identity of $S$.
One can check that $Sym(S)$ is the free commutative
monoid generated by $S\{-1\}$ in $\sseq A$. In other words,
for any commutative monoid object $M$ in the category of
symmetric sequences, any morphism from $S$ to $M_1$ in $\A$
can be extended uniquely into a morphism of commutative monoids
from $Sym(S)$ to $M$ in $\sseq A$ (exercise left to the reader).\\
\indent Given an object $S$ of $\A$, we define the
category of \emph{symmetric $S$-spectra} as the category of
modules in $\sseq A$ over the commutative monoid $Sym(S)$.
In explicit terms, a symmetric $S$-spectra $E$ is a collection $(E_{n},\s_{n})_{n\geq 0}$
where for any $n\geq 0$, $E_{n}$ is an object of $\A$
and $\s_{n}:S\otimes E_{n}\To E_{n+1}$
is a morphism in $\A$ such that the induced map obtained by composition
$$S^{\otimes m}\otimes E_{n}\To S^{\otimes m-1}\otimes E_{n+1}\To\dots\To
S\otimes E_{m+n-1}\To E_{m+n}$$
is $\mathfrak S_{m}\times\mathfrak S_{n}$-equivariant
 for any integers $n \geq 0, m \geq 0$.
A morphism of symmetric $S$-spectra $u:(E_{n},\s_{n})\To(F_{n},\tau_{n})$
 is a collection of $\mathfrak S_{n}$-equivariant maps
$u_{n}:E_{n}\To F_{n}$ such that the square
$$\xymatrix{
S\otimes E_{n}\ar[r]^{\s_{n}}\ar[d]_{S\otimes u_{n}}&E_{n+1}\ar[d]^{u_{n+1}}\\
S\otimes F_{n}\ar[r]_{\tau_{n}}&F_{n+1}
}$$
commutes in $\A$ for any integer $n \geq 0$.
\end{paragr}

\begin{paragr}\label{hypothesedebase}
We consider now a closed symmetric monoidal
Grothendieck abelian category $\A$
endowed with a weakly flat descent structure
$(\G,\H)$ and a flat set $\T$ of complexes of $\A$ (see \ref{defflatgcofcomploc}).
By virtue of Corollary \ref{cmf4}, the category $\Comp{(\A)}$
is endowed with a closed model category structure associated to $\T$
which is compatible with the tensor product. Let $S$
be a $\G$-cofibrant complex. We want $S$
to be invertible in the following sense: we want
the derived tensor product
by $S$ to be an equivalence of categories. For this purpose,
we will embed the category $\Comp{(\A)}$ in the category of
symmetric $S$-spectra by a symmetric monoidal functor,
and we will define a model category structure on the
category of symmetric $S$-spectra for which the tensor
product by $S$ will be a left Quillen equivalence.
We will have to study first the homotopy
theory of complexes of symmetric sequences associated to $\T$.\\
\indent The category $\Comp{(\A)}^{\mathfrak S}$
of symmetric sequences in $\Comp{(\A)}$ is canonically
equivalent to the category $\Comp{(\A^{\mathfrak S})}$
of complexes of $\A^{\mathfrak S}$. Define 
$\G^{\mathfrak S}$ to be the class of objects
$E\{-i\}$ for $E$ in $\G$ and $i\geq 0$, and $\H^{\mathfrak S}$
(resp. $\T^{\mathfrak S}$)
the class of complexes $H\{-i\}$ for $H$ in $\H$ (resp. $T\{-i\}$ for $T$ in $\T$) 
 and any integer $i\geq 0$.
\end{paragr}

\begin{lm}\label{gensymm}
$\G^{\mathfrak S}$ is a generating family of $\A^{\mathfrak S}$.
If moreover $\G$ is flat, then $\G^{\mathfrak S}$ is flat.
\end{lm}

\begin{proof}
We first have to show that $\G^{\mathfrak S}$ is a generating family
for $\A^{\mathfrak S}$. But for a symmetric sequence
$A_{*}$, we have
$$\Hom_{\A^{\mathfrak S}}(E\{-i\},A_{*})=\Hom_{\A}(E,A_{i}) \ .$$
As the family of evaluation functors $\sseq A \xrightarrow{Ev_n} \A$
is conservative and as $\G$ is a generating family of $\A$, this implies
that $\G^{\mathfrak S}$ is a generating family of $\A^{\mathfrak S}$.
It is easy to see that $\G^{\mathfrak S}$ is flat: this comes from the fact that for
any objects $X$ and $Y$ of $\A$ and any integers $i,j\geq 0$, we have
$X\{-i\}\otimes^{\mathfrak S}Y\{-j\}=(X\otimes Y)\{-i-j\}$.
\end{proof}

\begin{lm}\label{gcofsymmacycl}
If $C$ is a $\G$-cofibrant complex (resp. an acyclic complex)
then $C\{-i\}$ is $\G^{\mathfrak S}$-cofibrant (resp. acyclic).
\end{lm}

\begin{proof}
This follows from the fact that
the functors $\A\To\A^{\mathfrak S}$, $X\longmapsto X\{-i\}$
preserve colimits and are exact.
\end{proof}

\begin{lm}\label{symmadjoncelem}
For any complex
$C$ of $\A$ and any complex $K_{*}$ of $\A^{\mathfrak S}$, we have
canonical isomorphisms for $i\geq 0$
$$\Hom_{\Htp(\A^{\mathfrak S})}(C\{-i\},K_{*})\simeq\Hom_{\Htp(\A)}(C,K_{i}) \ .$$
\end{lm}

\begin{proof}
We have $\Hom_{\Comp(\A^{\mathfrak S})}(C\{-i\},K_{*})\simeq\Hom_{\Comp(\A)}(C,K_{i})$
and this isomorphism is compatible with the cochain homotopy relation.
\end{proof}

\begin{lm}\label{symmadjonc}
For any complex
$C$ of $\A$ and any complex $K_{*}$ of $\A^{\mathfrak S}$, we have
canonical isomorphisms for $i\geq 0$
$$\Hom_{\Der(\A^{\mathfrak S})}(C\{-i\},K_{*})\simeq\Hom_{\Der(\A)}(C,K_{i}) \ .$$
\end{lm}

\begin{proof}
This follows from these two facts:
the functor $\Comp(\A)\To\Comp(\A^{\mathfrak S})$, $C\longmapsto C\{-i\}$
is a left Quillen functor for the injective model structures on
$\Comp(\A)$ and $\Comp(\A^{\mathfrak S})$ and the functor
$\Comp(\A^{\mathfrak S})\To\Comp(\A)$,
$C_{*}\longmapsto C_{i}$
preserves quasi-isomorphisms. Indeed, this shows that it is sufficient
to prove the result for $K_{*}$ fibrant with respect to the injective
model structure. But then this reduces to Lemma \ref{symmadjoncelem}.
\end{proof}

\begin{prop}\label{symmflatdescent}
$(\G^{\mathfrak S},\H^{\mathfrak S})$ is a weakly flat descent structure on $\A^{\mathfrak S}$
and the family $\T^{\mathfrak S}$ of the $T\{-i\}$'s for $T$ in $\T$ and $i\geq 0$
is flat. Moreover a complex of symmetric sequences $K_{*}$ is $\G^{\mathfrak S}$-local
(resp. $\H^{\mathfrak S}$-flasque, resp. $\T^{\mathfrak S}$-local)
if and only if for any integer $i\geq 0$, $K_{i}$ is $\G$-local
(resp. $\H$-flasque, resp. $\T$-local).

If $\G$ is flat, then then descent structure $(\G^{\mathfrak S},\H^{\mathfrak S})$ is flat.
\end{prop}

\begin{proof}
It follows from Lemmata \ref{symmadjoncelem} and \ref{symmadjonc}
that a complex of symmetric sequences $K_{*}$ is $\G^{\mathfrak S}$-local
(resp. $\H^{\mathfrak S}$-flasque, resp. $\T^{\mathfrak S}$-local)
if and only if for any integer $i\geq 0$, $K_{i}$ is $\G$-local
(resp. $\H$-flasque, resp. $\T$-local).
This implies that $(\G^{\mathfrak S},\H^{\mathfrak S})$ is a descent structure.
We deduce easily from Lemma \ref{gcofsymmacycl} that it is weakly flat.
If moreover $\G$ is flat, then $(\G^{\mathfrak S},\H^{\mathfrak S})$ is a flat descent structure:
this comes from Lemmata \ref{gensymm} and \ref{gcofsymmacycl}
and from Proposition \ref{cmf333}.
\end{proof}

\begin{prop}
A morphism of complexes of symmetric sequences $u_{*}:K_{*}\To L_{*}$
is a $\T^{\mathfrak S}$-equivalence (resp. a fibration with respect to the
$\G^{\mathfrak S}$-model structure associated to $\T^{\mathfrak S}$)
if and only if for any integer
$i\geq 0$, the map $u_{i}:K_{i}\To L_{i}$ is a $\T$-equivalence
(resp. a fibration with respect to the
$\G$-model structure associated to $\T$).
\end{prop}

\begin{proof}
The characterization of fibrations follows from Proposition \ref{symmflatdescent}
and from Corollary \ref{characfibrations3}.
If $i\geq 0$ is an integer, the evaluation functor at $i$
$$\Comp(\A)^{\mathfrak S}\To\Comp(\A)\quad , \qquad A_{*}\longmapsto A_{i}$$
sends the $\G^{\mathfrak S}$-cofibrations (resp.
the trivial $\G^{\mathfrak S}$-cofibrations) to $\G$-cofibrations
(resp. trivial $\G$-cofibrations): the class of $\G^{\mathfrak S}$-cofibrations (resp.
of trivial $\G^{\mathfrak S}$-cofibrations) is generated by maps
of the form $u\{-j\}:A\{-j\} \To B\{-j\}$ where $j\geq 0$ is an integer and
$u$ is a $\G$-cofibration (resp. a trivial $\G$-cofibration);
see and apply \ref{defJcof} and Remark \ref{generators} to $\A^{\mathfrak S}$.
As the evaluation functor at $i$ preserves colimits, it is sufficient to
check that for any $\G$-cofibration (resp. trivial $\G$-cofibration)
$u:A\To B$ and any integer $j\leq 0$, the morphism
$u\{-j\}_{i}:A\{-j\}_{i}\To B\{-j\}_{i}$ is a $\G$-cofibration (resp. a trivial $\G$-cofibration).
But for $i\neq j$, everything is zero, so that we only have to
deal with the case $i=j$. Then the map $u\{-i\}_{i}$ is just
$$\bigoplus_{\mathfrak S_i}A=\mathfrak S_i \times A\To
\mathfrak S_i\times B=\bigoplus_{\mathfrak S_i}B \ . $$
Our claim thus follows from the fact that $\G$-cofibrations (resp. trivial $\G$-cofibrations)
are stable by (finite) direct sums. Let $u_{*}:K_{*}\To L_{*}$
be a morphism of symmetric sequences. Using the existence
of factorizations into a trivial cofibration followed by
a fibration in the $\G^{\mathfrak S}$-model structure associated to $\T^{\mathfrak S}$,
one can produce a diagram
$$\xymatrix{
{K_{*}}\ar[r]^k\ar[d]_{u}&M_{*}\ar[d]^v\\
{L_{*}}\ar[r]_{l}&N_{*}
}$$
in which $k$ and $l$ are trivial cofibrations and $M_{*}$
and $N_{*}$ are fibrant. As $k$ and $l$ are both termwise
$\T$-equivalences and $\T ^{\mathfrak S}$-equivalences,
$u$ is a termwise $\T$-equivalence (resp. a $\T ^{\mathfrak S}$-equivalence)
if and only if $v$ has the same property.
But as $M_{*}$ and $L_{*}$ are fibrant, $v$
is a $\T ^{\mathfrak S}$-equivalence if and only if it is
a quasi-isomorphism, and as the $M_{i}$'s and the $N_{i}$'s
are fibrant, $v$ is a termwise $\T$-equivalence if and only if
its is a termwise quasi-isomorphism. But we know that the evaluation functors
are exact and form a conservative family of functors, so that
we get that $v$ is a termwise quasi-isomorphism if and only if
it is a quasi-isomorphism. Hence we deduce that
$u$ is a $\T^{\mathfrak S}$-equivalence if and only if it is a termwise $\T$-equivalence.
\end{proof}

\begin{paragr}\label{defprojcmfspectra}
We can now consider a $\G$-cofibrant complex $S$ of $\A$.
We have a commutative monoid $Sym(S)$ in the category
of symmetric sequences of $\Comp(\A)$ (see \ref{defunivsymmmnoid}).
We will write $\Sp S(\A)$ for the category of symmetric $S$-spectra
(\emph{i.e.} of $Sym(S)$-modules) in $\Comp(\A)^{\mathfrak S}$; see
\ref{defunivsymmmnoid} again. The category of symmetric $S$-spectra
is a closed symmetric monoidal category (as any
category of modules over a commutative monoid
in a symmetric monoidal category; see e.g. \cite{MLa}).
We will write $\otimes^{}_{Sym(S)}$ or simply
$\otimes$ for the tensor product in $\Sp S(\A)$,
and $\unit^{}_{S}$ for the unit (\emph{i.e.} $\unit^{}_{S}=Sym(S)$
as a symmetric sequence).
We have a symmetric monoidal functor
$$\Comp(\A)^{\mathfrak S}\To\Sp S(\A)\quad ,
\qquad A_{*}\longmapsto Sym(S)\otimes^{\mathfrak S}A_{*}$$
which has a right adjoint that consists to forget the action of $Sym(S)$
$$U : \Sp S(\A)\To\Comp(\A)^{\mathfrak S}\quad ,
\qquad E\longmapsto E_{*}\ .$$
If we think of symmetric $S$-spectra as collections
$E=(E_{n},\s_{n})_{n\geq 0}$, then $U(E)$ is the symmetric sequence
$E_{*}=(E_{n})_{n\geq 0}$.

A morphism of symmetric $S$-spectra will be said to be a
\emph{$\T$-equivalence} (resp. a \emph{$\T$-fibration})
if its image by $U$ is a termwise $\T$-equivalence (resp. a termwise
fibration with respect to the $\G$-local model structure
associated to $\T$ on $\Comp(\A)$\footnote{The notion of $\T$-fibration
depends on $\T$ and $\G$, but we will neglect this as we will
always work with a fixed set $\G$ of generators.}).
\end{paragr}

\begin{prop}\label{projcmfspectra}
The category $\Sp S(\A)$ is a proper cellular
symmetric monoidal model category with the $\T$-equivalences (resp. the $\T$-fibrations)
as weak equivalences (resp. as fibrations). Furthermore, if $\G$ is flat,
then the monoid axiom holds in $\Sp S(\A)$.
\end{prop}

\begin{proof}
The first assertion is proved by Hovey~\cite[Theorem 8.2 and Theorem 8.3]{Hov3}.
The right properness follows from the
fact that the $\G^{\mathfrak S}$-model structure associated
to $\T^{\mathfrak S}$ on $\Comp(\A)^{\mathfrak S}$ is
proper and the fact that the functor $U$ preserves and detect
weak equivalences, fibrations and pullbacks.
We get the left properness in a similar way: this follows from
Proposition \ref{carrehtpcartTlc} applied to $\A^{\mathfrak S}$ and from
the facts that the cofibrations of $\Sp S(\A)$ are monomorphisms and
that the functor $U$ preserves and detects
weak equivalences, monomorphisms and pushouts.
The last assertion is an application of \cite[Theorem 4.1]{SS}.
\end{proof}

\begin{paragr}
We denote by $\Der_{\T}(\A,S)$ the localisation of the category
of symmetric $S$-spectra by the termwise $\T$-equivalences.
\end{paragr}

\begin{prop}\label{injcmfspectra}
The category $\Sp S(\A)$ is a proper cellular
model category with the $\T$-equivalences (resp. the monomorphisms)
as weak equivalences (resp. as cofibrations).
\end{prop}

\begin{proof}
It follows from Propositions \ref{carrehtpcartTlc} and
\ref{carrehtpcartTlc2} that the class of monomorphisms
which are $\T$-equivalences is stable by pushout, transfinite
composition and retract. As the class of $\T$-equivalences
is the class of weak equivalences of a cofibrantly generated model category,
it is accessible. We thus obtain the expected model
category structure by a straight application of Jeff Smith Theorem; see \cite{Beke}.
The left properness comes from the fact that any object is cofibrant for this
model structure. The right properness comes from the fact
that any fibration for this model structure is a $\T$-fibration
and from the right properness of the model structure of Proposition \ref{projcmfspectra}.
\end{proof}

\begin{prop}\label{niceSbasechange}
If $\G$ is flat, then the functor $A_{*}\longmapsto Sym(S)\otimes^{\mathfrak S}A_{*}$
preserves termwise $\T$-equivalences.
\end{prop}

\begin{proof}
As $S$ is $\G$-cofibrant, for any integer $n\geq 0$,
$S^{\otimes n}$ is $\G$-cofibrant as well. This implies that
the functor $A\longmapsto S^{\otimes n}\otimes A$
preserves $\T$-equivalences in $\Comp(\A)$:
to see this, just remember that the canonical map
$S^{\otimes n}\otimes^\derL A\To S^{\otimes n}\otimes A$
is an isomorphism in $\Der_{\T}(\A)$ (see \ref{tensdercoincide} and
the last assertion of \ref{cmf4}). This proposition thus follows from
the construction of $Sym(S)\otimes^{\mathfrak S}A_{*}$
(see \ref{defsymtensoreq}) and from the fact that $\T$-equivalences are
stable by finite direct sums.
\end{proof}

\begin{rem}\label{termwisepectratriang}
The model categories of Propositions \ref{projcmfspectra} and \ref{injcmfspectra}
are stable. This means that a commutative square in the category of symmetric
$S$-spectra is a homotopy pushout square if and only if it is a homotopy
pullback square. In particular, the homotopy category $\Der_{\T}(\A,S)$
(localization of $\Sp S(\A)$ by the termwise $\T$-equivalences)
is canonically endowed with a symmetric monoidal triangulated
category structure, and the functor
$$\Der_{\T^{\mathfrak S}}(\A^{\mathfrak S})\To \Der_{\T}(\A,S)\quad ,
\qquad A_{*}\longmapsto Sym(S)\otimes^{\mathfrak S}A_{*}$$
is a symmetric monoidal functor.
\end{rem}

\begin{paragr}
We have a symmetric monoidal functor
$$\Sigma^{\infty} : \Comp(\A)\To\Sp S(\A)\quad , \qquad
A\longmapsto Sym(S)\otimes^{\mathfrak S} A\{0\} \ .$$
In particular, we have the formula
\begin{equation}\label{sigmainftymonoidal}
\Sigma^\infty(A\otimes B)\simeq\Sigma^{\infty}(A)\otimes\Sigma^{\infty}(B) \ .
\end{equation}
The functor $\Sigma^\infty$ is a left adjoint of the functor
$$\Omega^\infty : \Sp S(\A)\To\Comp(\A)\quad , \qquad
(E_{n},\s_{n})_{n\geq 0}\longmapsto E_{0} \ .$$
These two functors form a Quillen adjunction
if we consider the $\G$-model structure associated to $\T$ on $\Comp(\A)$
and the model structure of Proposition \ref{projcmfspectra} on $\Sp S(\A)$.
We introduce now a very ugly abuse of notation.
For an object $A$ of $\Comp(\A)$ and a symmetric $S$-spectrum $E$,
we define $E\otimes A$ to be $E\otimes\Sigma^\infty(A)$.
This can be rewritten as
\begin{equation}\label{Ssymmaction}
E\otimes A = E\otimes^{}_{Sym(S)}\Big(Sym(S)\otimes^{\mathfrak S}A\{0\}\Big)
\simeq E\otimes^{\mathfrak S}A\{0\} \ .
\end{equation}
In terms of collections $E=(E_{n},\s_{n})_{n\geq 0}$, we have
$E\otimes A=(E_{n}\otimes A,\s_{n}\otimes 1_{A})_{n\geq 0}$.
This makes our purpose more precise: we would like the functor
$E\longmapsto E\otimes S$ to be a left Quillen equivalence.
This is not the case in general, so that we have to localize a little further.

We will say that a symmetric $S$-spectrum $E=(E_{n},\s_{n})_{n\geq 0}$
is an \emph{$\Omega^\infty$-spectrum} if for any integer $n\geq 0$,
$E_{n}$ is a $\G$-local and $\T$-local complex of $\A$ and the map
induced by adjunction from $\s_{n}$
\begin{equation}\label{omegamaps}
\tilde{\s}_{n}:E_{n}\To\sHom(S,E_{n+1})
\end{equation}
is a quasi-isomorphism (here $\sHom$
is the internal Hom of the category of complexes of $\A$).
We remark that as $S$ is $\G$-cofibrant
and as $E_{n+1}$ is a fibrant object
of the $\G$-model structure associated to $\T$,
$\sHom(S,E_{n+1})$ is also fibrant. This implies that
the map \ref{omegamaps} is a quasi-isomorphism if
and only if it is a $\T$-equivalence.

A morphism of symmetric $S$-spectra $u : A\To B$
is a \emph{stable $\T$-equivalence} if for any
$\Omega^\infty$-spectrum $E$, the induced map
$$u^* : \Hom_{\Der_{\T}(\A,S)}(B,E)\To\Hom_{\Der_{\T}(\A,S)}(A,E)$$
is an isomorphism of abelian groups. For example, any termwise $\T$-equivalence
is a stable $\T$-equivalence.

A morphism of symmetric $S$-spectra is a \emph{stable $\G$-cofibration}
is it is a cofibration of the model structure of Proposition \ref{projcmfspectra}.

A morphism of symmetric $S$-spectra is a \emph{stable $\T$-fibration}
if it has the right lifting property with
respect to the stable $\G$-cofibrations which are also stable
$\T$-equivalences.
\end{paragr}

\begin{prop}\label{stablecmf}
The category of symmetric $S$-spectra is a stable proper cellular symmetric
monoidal model category with the stable $\T$-equivalences as weak equivalences,
the stable $\G$-cofibrations as cofibrations, and the stable
$\T$-fibrations as fibrations. Moreover, the stable $\T$-fibrations
are the termwise $\G$-surjective morphisms $p:E\To F$
such that $\mathrm{ker}\, p$ is an $\Omega^\infty$-spectrum.
In particular, the fibrant objects of this model structure
are the $\Omega^\infty$-spectra.
\end{prop}

\begin{proof}
It follows from \cite[Theorem 8.11]{Hov3} that
we have defined a model category structure.
As this model structure is by definition a left Bousfield
localization of a left proper cellular model category,
it is left proper. This model structure
is stable as the stable $\T$-equivalences are the maps
inducing isomorphisms on the cohomological
functors $\Hom_{\Der_{\T}(\A,S)}(?,E)$ with $E$
any $\Omega^\infty$-spectrum. The right properness
follows then from the right properness of the
model category of Proposition \ref{projcmfspectra}
as follows. Consider a pullback square of symmetric
spectra
$$\xymatrix{
E\ar[r]^{u}\ar[d]_{p}&F\ar[d]^q\\
E'\ar[r]_{u'}&F'
}$$
where $p$ is a stable $\T$-fibration.
As $p$ is also a fibration for the model structure of Proposition \ref{projcmfspectra},
this square is a homotopy pullback square. But as this latter
model structure is stable, this square is also a homotopy
pushout square. This implies that this square is a homotopy
pushout square in the model category defined by the
stable $\T$-equivalences. Again, as this latter model structure is
stable, we get that this commutative square has to be
a homotopy pullback square. Hence, if moreover $u'$
is a stable $\T$-equivalence, $u$ has to be a stable $\T$-equivalence
as well, and we have proved the right properness.
The characterization of the stable $\T$-fibrations
starts with the description of the fibrant objects as the
$\Omega^\infty$-spectra which comes from \cite[Theorem 8.8]{Hov3}.
We can then finish the proof using the fact that the fibrations of the model structure
of Proposition \ref{projcmfspectra} are exactly
the termwise $\G$-surjective maps with terwise $\G$-local
and $\T$-local kernels (see \ref{characfibrations2}) and
reproducing the proof of Proposition \ref{fibrationslocalisees}.
\end{proof}

\begin{paragr}\label{defSstablecmf}
The model structure of Proposition \ref{stablecmf} will be called
the \emph{$S$-stable model structure associated to $\T$}.
We denote by $\Der_{\T}(\A,S^{-1})$ the localization of the
category $\Sp S(\A)$ by the stable $\T$-equivalences.
\end{paragr}

\begin{paragr}\label{defSmoinsun}
It remains to show that the functor
$E\longmapsto E\otimes S=E\otimes\Sigma^\infty(S)$
is a left Quillen equivalence. We already know that
it is a left Quillen functor with respect to the $S$-stable model
category structure associated to $\T$. We write
$$\Der_{\T}(\A,S^{-1})\To\Der_{\T}(\A,S^{-1})\quad , \qquad
E\longmapsto E\otimes^\derL S$$
for the corresponding total left derived functor. We define
the symmetric $S$-spectra
$S^{-1}=Sym(S)\otimes^{\mathfrak S}S\{-1\}$. The tensor product
by $S^{-1}$ 
$$\Sp S(\A)\To\Sp S(\A)\quad , \qquad E\longmapsto E\otimes S^{-1}
\simeq E\otimes^{\mathfrak S} S\{-1\}$$
is also a left Quillen functor and has also a total left derived functor
$$\Der_{\T}(\A,S^{-1})\To\Der_{\T}(\A,S^{-1})\quad , \qquad
E\longmapsto E\otimes^\derL S^{-1} \ .$$
Note that, if moreover $\G$ is flat,
as both $\Sigma^\infty(S)$ and $S^{-1}$ are stably $\G$-cofibrant, for any
symmetric $S$-spectrum $E$, we have canonical
isomorphisms in $\Der_{\T}(\A,S^{-1})$ (\ref{stablemonoidaxiom02})
$$E\otimes^\derL S\simeq E\otimes S\quad\text{and}\quad
E\otimes^\derL S^{-1}\simeq E\otimes S^{-1} $$
(but this property will not be needed).
\end{paragr}

\begin{prop}\label{stabiliteenfin}
The functors $E\longmapsto E\otimes^\derL S$
and $E\longmapsto E\otimes^\derL S^{-1}$ are equivalence of
categories and are quasi-inverse to each other.
In particular, the functor $E\longmapsto E\otimes S=E\otimes\Sigma^\infty(S)$
is a left Quillen equivalence.
\end{prop}

\begin{proof}
See \cite[Theorem 8.10]{Hov3}.
\end{proof}

\begin{ex} \label{DM}
Consider the notations of \ref{DMeff}.

Let $\GG_m$ be the multiplicative group scheme over $S$.
It is a pointed scheme with $S$-point the unit section $s$.
We associate to this pointed scheme a sheaf with transfers 
$\tilde L_S(\GG_m)$ equals to the cokernel of the split monomorphism
$L_S(S) \xrightarrow{s_*} L_S(\GG_m)$.

A symmetric $\tilde L_S(\GG_m)$-spectrum in the category 
$\Comp(\A)^{\mathfrak S}$ (cf \ref{defprojcmfspectra})
will be called a motivic spectrum over $S$.

According to Proposition \ref{stablecmf},
the category of symmetric Tate spectra $\spt S$ over $S$
is a stable symmetric monoidal category. We denote
its homotopy category by $\DM(S)$.

By construction, we obtain a pair of adjoint functors
$$
\Sigma^{\infty}:\DMe(S)
 \leftrightarrows \DM(S):\Omega^\infty
$$
such that $\Sigma^\infty$ is symmetric monoidal,
and $\Sigma^\infty \tilde M_S(\GG_m)$ (usually denoted by $\ZZ_S(1)[1]$)
is invertible
(cf \ref{stabiliteenfin})\footnote{In fact,
 we can state a universal property satisfied by 
 $\DM(S)$ from these two properties}.

Using an obvious extension of the base change 
(resp. forget the base) functor from complexes to spectra,
and applying the obvious analog of Proposition \ref{fonct2}, 
we can construct for any morphism of regular 
schemes $f:T \rightarrow S$ a pair of adjoint functors
$$
\derL f^*:\DM(S) \leftrightarrows \DM(T):\derR f_*
$$
and in the case where $f$ is smooth another pair
$$
\derL f_\sharp:\DM(S) \leftrightarrows \DM(T):\derR f^*=f^*.
$$

This constructions will be generalized and studied more closely in the forthcoming paper
\cite{DMCD}.
\end{ex}

\begin{paragr}\label{suspdecalage}
We will need a little more details about the relationship
between the tensor product by $S$ and $S^{-1}$.
We will write $\sHom_{\mathfrak S}$ and $\sHom_{S}$
for the internal Hom's in the categories
of symmetric sequences and of symmetric $S$-spectra
respectively.

Let $E=(E_{n},\s_{n})_{n\geq 0}$ be a symmetric $S$-spectrum.
One can define the symmetric $S$-spectrum $E\{1\}$ by
$$E\{1\}=\sHom_{\mathfrak S}(\unit\{-1\},E)\simeq
\sHom_{S}(Sym(S)\{-1\},E) \ .$$
In other words, $E\{1\}=(E\{1\}_{n},\s\{1\}_{n})_{n\geq 0}$, where
for $n\geq 0$, we have $E\{1\}_{n}=E_{n+1}$ and
$\s\{1\}_{n}=\s_{n+1}$, and the action of $\mathfrak S_{n}$
on $E_{n+1}$ is induced by the obvious inclusion
$\mathfrak S_{n}\subset\mathfrak S_{n+1}$. The maps
$\s_{n}$ define a morphism of symmetric $S$-spectra
called the \emph{suspension map}
\begin{equation}\label{suspdec1}
\s : E\otimes S\To E\{1\} \ .
\end{equation}
By adjunction, this defines a map
\begin{equation}\label{suspdec2}
\tilde{\s}:E\To\sHom_{S}(\Sigma^\infty(S),E\{1\}) \ .
\end{equation}
But $\sHom_{S}(\Sigma^\infty(S),E\{1\})\simeq
\sHom_{S}(Sym(S)\otimes^{\mathfrak S} S\{-1\},E)$, so that
we have a map
\begin{equation}\label{suspdec3}
\tilde{\s}:E\To\sHom_{S}(Sym(S)\otimes^{\mathfrak S} S\{-1\},E) \ .
\end{equation}
Then $E$ is an $\Omega^\infty$-spectrum if and only if
it is termwise $\G$-local and $\T$-local and if the map \eqref{suspdec3}
is a termwise $\T$-equivalence (this is an easy translation of
the definition). We deduce from this
and from the fact that the functor $?\otimes S$ preserves
stable $\T$-equivalences (see \ref{stablemonoidaxiom02})
that for any $\Omega^\infty$-spectrum $E$,
the map \eqref{suspdec1} is a stable $\T$-equivalence\footnote{The proof
of Proposition \ref{stabiliteenfin} follows from this fact once we have noticed that
$\sHom_{S}(Sym(S)\otimes^{\mathfrak S} S\{-1\},?)$ is a right adjoint
to the functor $?\otimes S^{-1}$.}. As $E\{1\}$
is then still an $\Omega^\infty$-spectrum, we can iterate
this process and obtain for any $\Omega^\infty$-spectrum $E$
a canonical stable $\T$-equivalence
\begin{equation}\label{suspdec4}
\s^n : E\otimes S^{\otimes n}=E\otimes\Sigma^\infty(S)^{\otimes n}\To E\{n\} \ .
\end{equation}
\end{paragr}

\begin{paragr}\label{suspensioninfinie}
The functor $\Sigma^\infty:\Comp(\A)\To\Sp S(\A)$
is a left Quillen functor from the $\G$-model structure associated to
$\T$ to the $S$-stable model structure associated to $\T$ (this is
immediate).
If moreover $\G$ is flat, it follows easily from Proposition \ref{niceSbasechange}
that it sends $\T$-equivalences to stable $\T$-equivalences.

The fact $\Sigma^\infty$ is a left quillen functor
implies that his right adjoint $\Omega^\infty$
is a right Quillen functor and that the total left derived functor
$$\derL\Sigma^\infty : \Der_{\T}(\A)\To\Der_{\T}(\A,S^{-1})$$
is a left adjoint of the total right derived functor of $\Omega^\infty$
$$\derR\Omega^\infty : \Der_{\T}(\A,S^{-1})\To\Der_{\T}(\A) \ .$$
In particular, for any complex $X$ of $\A$ and any symmetric $S$-spectrum $E$,
we have canonical isomorphisms of abelian groups
$$\Hom_{\Der_{\T}(\A,S^{-1})}(\derL\Sigma^\infty(X),E)\simeq
\Hom_{\Der_{\T}(\A)}(X,\derR\Omega^\infty(E)) \ .$$
Moreover, for any $\Omega^\infty$-spectrum $E$, we have
$R\Omega^\infty(E)=E_{0}$.

A \emph{weak $\Omega^\infty$-spectrum} is a symmetric $S$-spectrum
$E=(E_{n},\s_{n})_{n\geq 0}$ such that for any $n\geq 0$, the map
$$E_{n}\To\derR\sHom(S,E_{n+1})$$
obtained by adjunction in $\Der_{\T}(\A,S)$
from the map $\sigma_{n}:E_{n}\otimes S\To E_{n+1}$ is an isomorphism.
It is obvious that any $\Omega^\infty$-spectrum
is a weak $\Omega^\infty$-spectrum.
\end{paragr}

\begin{prop}\label{weakomega}
Let $E$ be a weak $\Omega^\infty$-spectrum and $n\geq 0$ an integer.
Then we have canonical
isomorphisms $\derR\Omega^\infty(E\otimes S^{\otimes n})\simeq E_{n}$ in
$\Der_{\T}(\A)$. In other words, for any object $X$ of $\Der_{\T}(\A)$,
we have canonical isomorphisms
$$\begin{aligned}
\Hom_{\Der_{\T}(\A)}(X,E_{n})
&\simeq\Hom_{\Der_{\T}(\A)}(X,\derR\Omega^\infty(E\otimes S^{\otimes n}))\\
&\simeq\Hom_{\Der_{\T}(\A,S^{-1})}(\derL\Sigma^\infty(X),E\otimes S^{\otimes n}) \ .
\end{aligned}
$$
Moreover, the iterated suspension map
$\s^n : E\otimes S^{\otimes n}\To E\{n\}$
is an isomorphism in $\Der_{\T}(\A,S^{-1})$.
\end{prop}

\begin{proof}
We know the assertion is true when $E$ is actually an
$\Omega^\infty$-spectrum (see \ref{suspdecalage}).
In general, thanks to the
existence of the model structure of
Proposition \ref{projcmfspectra}, one can choose a
termwise $\T$-equivalence $E\To E'$ such that
$E'$ is termwise $\G$-local and $\T$-local.
But $E'$ is also a weak $\Omega^\infty$-spectrum,
and this implies immediately that $E'$ is an $\Omega^\infty$-spectrum.
As the maps $E_{i}\To E'_{i}$ are $\T$-equivalences for any $i\geq 0$, this
implies that $E\{n\}\To E'\{n\}$ is also a termwise $\T$-equivalence
for any integer $n$.
Hence we are reduced to the $\Omega^\infty$-spectrum case.
\end{proof}

\begin{prop}\label{carrehtpcartTlc3}
A pushout square of symmetric $S$-spectra
$$\xymatrix{
E\ar[r]^{u}\ar[d]_{i}&F\ar[d]^j\\
E'\ar[r]_{u'}&F'
}$$
in which $i$ is a monomorphism is a homotopy pushout square
with respect to the $S$-stable model structure associated to $\T$.
In particular, if $i$ (resp. $u$) is a stable $\T$-equivalence, then
$j$ (resp. $u'$) is a stable $\T$-equivalence.
\end{prop}

\begin{proof}
It is sufficient to prove that this commutative square is a homotopy pushout
square for the model structure of Proposition \ref{projcmfspectra}.
This square is obviously a homotopy pushout square for the model
structure of Proposition \ref{injcmfspectra}, which implies
that it has to be a homotopy pushout square for the model
structure of Proposition \ref{projcmfspectra} as well.
\end{proof}

\begin{prop}\label{carrehtpcartTlc4}
The stable $\T$-equivalences are stable by filtering
colimits in $\Sp S(\A)$.
\end{prop}

\begin{proof}
We know that the termwise $\T$-equivalences
are stable by filtering colimits (this follows immediately
from \ref{carrehtpcartTlc2}). This is equivalent to say that
for any filtering category $I$ and any functor $F$
from $I$ to $\Sp S(\A)$, the canonical map
$$\hocolim_{i}F_{i}\To\varinjlim_{i}F_{i}$$
is a termwise $\T$-equivalence, hence
a stable $\T$-equivalence. But as the functor
$\hocolim$ sends stable $\T$-equivalences to
stable $\T$-equivalences, this achieves the proof.
\end{proof}

\begin{prop}\label{injstablecmf}
The category of symmetric $S$-spectra
is endowed with a cellular proper model category
structure with the stable $\T$-equivalences as weak equivalences
and the monomorphisms as cofibrations.
\end{prop}

\begin{proof}
It follows from Propositions \ref{carrehtpcartTlc3} and
\ref{carrehtpcartTlc4} that the class of monomorphisms
which are stable $\T$-equivalences is stable by pushout, transfinite
composition and retract. We thus obtain the expected model structure
as a left Bousfield localization of the model structure
of Proposition \ref{injcmfspectra}. The left properness is obvious
(any object is cofibrant) and the right properness comes from
the right properness of the $S$-stable model structure
associated to $\T$ and from the fact that any fibration
is in particular a fibration for the $S$-stable model structure.
\end{proof}

\begin{paragr}
The model structure of Proposition \ref{injstablecmf} will be called the
\emph{injective $S$-stable model structure associated to $\T$}.

From now on, we suppose that $\G$ is flat.
\end{paragr}

\begin{lm}\label{stablemonoidaxiom00}
Let $A_{*}$ be a $\G^{\mathfrak S}$-cofibrant
symmetric sequence, and $E$ be a symmetric $S$-spectrum.
Then the canonical map
$$E\otimes^\derL\big(Sym(S)\otimes^{\mathfrak S} A_{*}\big)
\To E\otimes\big(Sym(S)\otimes^{\mathfrak S} A_{*}\big)$$
is an isomorphism in $\Der_{\T}(\A,S)$ (hence in
$\Der_{\T}(\A,S^{-1})$ as well).
\end{lm}

\begin{proof}
We begin with the following observation: the functor
$$\Comp(\A)^{\mathfrak S}\To\Comp(\A)^{\mathfrak S}
\quad , \qquad E\longmapsto E\otimes^{\mathfrak S} A_{*}$$
preserves termwise $\T$-equivalences. This is because
by virtue of Lemma \ref{gensymm},
$\G^{\mathfrak S}$ is a flat generating family of $\A^{\mathfrak S}$; see
Remark \ref{tensdercoincide}.

Let now $E$ be a symmetric $S$-spectrum.
We can choose a termwise $\T$-equivalence $E'\To E$ in
$\Sp S (\A)$ with $E'$ cofibrant for the model structure
of Proposition \ref{projcmfspectra}. But then
we have a termwise $\T$-equivalence (hence a stable
$\T$-equivalence)
$$E'\otimes\big(Sym(S)\otimes^{\mathfrak S} A_{*}\big)\simeq
E'\otimes^{\mathfrak S} A_{*}\To
E\otimes^{\mathfrak S} A_{*}\simeq
E\otimes\big(Sym(S)\otimes^{\mathfrak S} A_{*}\big) \ .$$
This implies our assertion as both $E'$ and $Sym(S)\otimes^{\mathfrak S} A_{*}$
are cofibrant symmetric $S$-spectra.
\end{proof}

\begin{lm}\label{stablemonoidaxiom01}
Let $E$ be a symmetric $S$-spectrum. Then the functor
$$\Sp S(\A)\To\Sp S(\A)\quad , \qquad
F\longmapsto E\otimes F$$
is a left Quillen functor from the model structure of Proposition \ref{projcmfspectra}
to the model structure of Proposition \ref{injcmfspectra}.
\end{lm}

\begin{proof}
 Proposition \ref{monoidaxiom} applied to $\A^{\mathfrak S}$
 gives the following.
If $u : A_{*}\To B_{*}$ is a $\G^{\mathfrak S}$-cofibration of symmetric sequences,
then the induced map
$$E\otimes\big(Sym(S)\otimes^{\mathfrak S} A_{*}\big)\simeq
E'\otimes^{\mathfrak S} A_{*}\To
E\otimes^{\mathfrak S} B_{*}\simeq
E\otimes\big(Sym(S)\otimes^{\mathfrak S} B_{*}\big)$$
is a monomorphism. If moreover $u$ is a quasi-isomorphism,
so is the induced map above. Hence the lemma is proved
for $\T=\varnothing$. To prove the general case, it is sufficient
to prove that the tensor product by $E$
sends the generating trivial cofibrations of the model
category of Proposition \ref{projcmfspectra}
to termwise $\T$-equivalences. But the generating
trivial cofibrations are all of shape
$$Sym(S)\otimes^{\mathfrak S}A\{-i\}\To
Sym(S)\otimes^{\mathfrak S}B\{-i\}$$
where $i\geq 0$ is an integer, and $A\To B$ is
a $\G$-cofibration and a $\T$-equivalence in $\Comp(\A)$.
We conclude easily using Lemma \ref{stablemonoidaxiom00}.
\end{proof}

\begin{prop}\label{stablemonoidaxiom02}
Let $E$ and $F$ be symmetric $S$-spectra.
If $E$ or $F$ is stably $\G$-cofibrant (\emph{i.e.} a
cofibrant object of the model category of Proposition \ref{projcmfspectra}),
then the canonical map
$$E\otimes^\derL F\To E\otimes F$$
is an isomorphism in $\Der_{\T}(\A,S)$ (hence in
$\Der_{\T}(\A,S^{-1})$ as well).
\end{prop}

\begin{proof}
Let $E$ be a fixed symmetric $S$-spectrum.
This object gives rise to two triangulated endofunctors
of the category $\Der_{\T}(\A,S)$.
We denote the first one by $\Phi$. It is defined
by $\Phi(F)=E\otimes^\derL F$ (where we insist that
$\otimes^\derL$ is the total left derived functor of the
left Quillen bifunctor $\otimes$ on $\Sp S(\A)$).
The second one, denoted by $\Psi$
 is the total left derived functor
of the functor $F\longmapsto E\otimes F$ (which is well
defined by virtue of Lemma \ref{stablemonoidaxiom01}).
We have a canonical map of triangulated functors
$\Phi\To\Psi$, and we would like it to be an isomorphism.
We already know by Lemma \ref{stablemonoidaxiom00}
that it induces an isomorphism
on objects of type $Sym(S)\otimes^{\mathfrak S}A_{*}$
for any $\G^{\mathfrak S}$-cofibrant symmetric sequence
$A_{*}$. But as the cokernels of the generating
cofibrations of the model structure of Proposition \ref{projcmfspectra}
are all of this type, it follows from \cite[Theorem 7.3.1]{Hov}
that the objects of shape $Sym(S)\otimes^{\mathfrak S}A_{*}$
form a generating family of the triangulated category $\Der_{\T}(\A,S)$.
This means that the smallest localizing
subcategory of $\Der_{\T}(\A,S)$ that contains the $Sym(S)\otimes^{\mathfrak S}A_{*}$'s
is $\Der_{\T}(\A,S)$ itself. This implies that $\Phi\To\Psi$
is an isomorphism: it just remains to
prove that the full sucategory of $\Der_{\T}(\A,S)$ that consists
of object $F$ such that $\Phi(F)\simeq\Psi(F)$ is a triangulated subcategory
and is stable by direct sums;
which is obvious. But if $F$ is stably $\G$-cofibrant, the map
$\Phi(F)\To\Psi(F)$ is just the map $E\otimes^\derL F\To E\otimes F$.
\end{proof}

\begin{prop}\label{stablemonoidaxiom0}
For any symmetric $S$-spectrum $E$, the functor
$A\longmapsto A\otimes E$ is a left Quillen functor from
the $S$-stable model structure associated to $\T$
to the injective $S$-stable model structure associated to $\T$.
\end{prop}

\begin{proof}
This follows immediately from Lemmata \ref{1locbousfieldcmf} and
\ref{stablemonoidaxiom01} and from Proposition
\ref{stablemonoidaxiom02}.
\end{proof}

\begin{cor}\label{stablemonoidaxiom}
If $\G$ is flat, then
the $S$-stable model category structure associated to $\T$
satisfies the monoid axiom of \cite{SS}.
\end{cor}

\begin{proof}
This is an immediate consequence of
Proposition \ref{stablemonoidaxiom0}.
\end{proof}

\begin{paragr}
We suppose that $\G$ is flat until the end of this section.

An \emph{$S$-ring spectrum}, or simply a \emph{ring spectrum}
is a commutative monoid in the category of symmetric $S$-spectra.
Such a ring spectrum is \emph{commutative} if it is
commutative as a monoid.

If $\ste$ is a ring spectrum, we write $\Sp S (\A,\ste)$
for the category of left $\ste$-modules. If moreover $\ste$ is commutative,
the category $\Sp S (\A,\ste)$ is a closed symmetric monoidal
category. We have a functor
$$\Sp S(\A)\To \Sp S (\A,\ste)\quad , \qquad
F\longmapsto \ste\otimes F$$
which is a left adjoint to the forgetful functor
$$\Sp S(\A,\ste)\To \Sp S (\A)\quad , \qquad
M\longmapsto M$$
If $\ste$ is commutative, the functor $F\longmapsto \ste\otimes F$
is of course symmetric monoidal. A morphism of $\ste$-modules
is a \emph{stable $\T$-equivalence} (resp. a \emph{stable $\T$-fibration})
if it is so as a morphism of symmetric $S$-spectra.

We will write $\Der_{\T}(\A,\ste)$ for the localization
of the category of $\ste$-modules by the class of stable $\T$-equivalences.
\end{paragr}

\begin{cor}\label{cmfstablemod}
Let $\ste$ be a ring spectrum. Then the category of left
$\ste$-modules in $\Sp S(\A)$ is a cellular proper stable 
model category with the stable $\T$-equivalences
as weak equivalences and the stable $\T$-fibrations
as fibrations. If moreover $\ste$ is commutative, this is
a symmetric monoidal model category which satisfies the
monoid axiom.
\end{cor}

\begin{proof}
The fact that we get a cofibrantly generated
model category (symmetric monoidal and satisfying
the monoid axiom if $\ste$ is commutative) comes from
\cite[Theorem 4.1]{SS}. The proof of properness
and of stability is left to the reader
(it is the same as in the proof of \ref{projcmfspectra}).
\end{proof}

\bibliographystyle{smfalpha}
\bibliography{common}

\end{document}